\documentclass[12pt]{article}

\usepackage[mathscr]{eucal}
\usepackage{amsmath, amsthm, amssymb, latexsym, graphicx,courier,makeidx,subeqnarray,color}
\usepackage{bm}
\usepackage{color}

\newtheorem{Theorem}{Theorem}

\newtheorem{Lemma}[Theorem]{Lemma}

\newtheorem{Problem}[Theorem]{Problem}

\newtheorem{Remark}{Remark}[section]

\def\real{\mathbb{R}}
\def\bar{\overline}

\newcommand{\bu}{\mbox{\boldmath{$u$}}}
\newcommand{\bv}{\mbox{\boldmath{$v$}}}
\newcommand{\bw}{\mbox{\boldmath{$w$}}}
\newcommand{\bx}{\mbox{\boldmath{$x$}}}

\newcommand{\fb}{\mbox{\boldmath{$f$}}}

\newcommand{\hb}{\mbox{\boldmath{$h$}}}
\newcommand{\bxi}{\mbox{\boldmath{$\xi$}}}

\newcommand{\bsigma}{\mbox{\boldmath{$\sigma$}}}
\newcommand{\btau}{\mbox{\boldmath{$\tau$}}}
\newcommand{\bvarepsilon}{\mbox{\boldmath{$\varepsilon$}}}
\newcommand{\bnu}{\mbox{\boldmath{$\nu$}}}
\newcommand{\bzeta}{\mbox{\boldmath{$\zeta$}}}

\newcommand{\bzero}{\mbox{\boldmath{$0$}}}

\def\lista#1
{{ \itemindent 0.0cm \labelsep .2cm \leftmargin 0.8cm \rightmargin
		0.0cm \labelwidth 0.6cm \topsep 0.0mm
		\parsep 0.0mm
		\itemsep 0.0mm
		\begin{list}{}
			{ \setlength{\leftmargin}{.8cm} \setlength{\rightmargin}{0.0cm}
				\setlength{\parsep}{0.0mm} \setlength{\topsep}{.0mm}
				\setlength{\parskip}{.0cm} \setlength{\itemsep}{.0cm} }
			{#1}\end{list}} }

\begin{document}
	
\title{\bf A Class of History-Dependent 
	Systems of Evolution Inclusions with Applications
	\thanks{\, 
	This project has received funding from the European Union's Horizon 2020 Research and Innovation Programme under the Marie Sk{\l}odowska-Curie grant agreement No. 823731 CONMECH. 
	It is supported by 
	Natural Science Foundation of Guangxi (Grant No: 2018GXNSFAA281353), 
	Beibu Gulf University Project No. 2018KYQD06, and the projects financed by the Ministry of Science and
	Higher Education of Republic of Poland under
	Grants Nos. 4004/GGPJII/H2020/2018/0 and
	440328/PnH2/2019, and the National Science Centre of Poland under Project No. 2021/41/B/ST1/01636.
}}

\author{
	Stanis{\l}aw Mig\'orski 
	\footnote{\,College of Applied Mathematics, Chengdu University of Information Technology, Chengdu 610225, Sichuan Province, P.R. China, and Jagiellonian University in Krakow, Chair of Optimization and Control, ul. Lojasiewicza 6, 30348 Krakow, Poland.
	Tel.: +48-12-6646666. E-mail address: stanislaw.migorski@uj.edu.pl.}
}

\renewcommand{\thefootnote}{\fnsymbol{footnote}}

\date{}
\maketitle
\thispagestyle{empty}

\begin{abstract}\noindent
This paper is devoted to studying a system of coupled nonlinear first order history-dependent evolution inclusions in the framework of 
evolution triples of spaces.
The multivalued terms are of the Clarke subgradient or of the convex subdifferential form.  
Using a surjectivity result for multivalued maps and a fixed point argument for a history-dependent operator, we prove that the system has a unique solution. 
We conclude with two examples of an evolutionary differential variational-hemivariational inequality and of a dynamic frictional contact problem in mechanics, which illustrate the abstract results.


\end{abstract}

\section{Introduction}\label{s1}
\setcounter{equation}0

In this paper we study the following system of nonlinear evolution inclusions:
\begin{eqnarray}
&&
w'(t) + A(t, \theta(t), w(t)) + 
({R}_1 w)(t) + 
\partial J (t, \theta(t), ({\mathcal S} w)(t), w(t)) 
\label{001}
\\[1mm]
&&
\qquad\qquad\qquad\qquad \quad \ \   
\, + \, \partial_c \varphi (t, \theta(t), 
({R} w)(t), w(t)) \ni h_1(t) 
\ \ \mbox{\rm a.e.} \ \ t \in (0,T), \nonumber \\[1mm]
&&
\theta'(t) + B(t, w(t), (R_2w)(t), \theta(t))
+ \partial g(t, w(t), \theta(t)) \ni h_2(t) \ \ \mbox{\rm a.e.} \ \ t \in (0,T), \label{002} \\[1mm]
&&
w(0) = w_0, \ \theta(0) = \theta_0, \label{003}
\end{eqnarray}
where $A$ and $B$ are nonlinear operators, $\partial J$ and $\partial g$ denote the Clarke subgradient with respect to the last variable of functions $J$ and $g$, respectively, $\partial_c \varphi$ stands for the subdifferential in the last variable of a convex function $\varphi$, and $R$, $R_1$,$R_2$ and $S$ are the so-called history-dependent operators.  
Our main goal is to prove existence and uniqueness of solution to the system as well as to provide applications 
to a new class of differential variational-hemivariational inequalities and to a new model of contact problems in mechanics.

System of the form (\ref{001})--(\ref{003}) appears naturally in the study of inequality problems in 
many areas, for instance, 
in mechanics where they represent the weak formulation of contact problems, see~\cite{HS,MOSBOOK,SST,SM2}. The convex function $\varphi$ and, in general, nonconvex functions $J$ and $g$ may have different meanings: they may represent 
various convex and nonconvex constitutive laws, as well as numerous contact and friction boundary conditions. 
Because of the presence of subdifferential terms, inclusion (\ref{001}) can be equivalently formulated as a variational-hemivariational inequality, and  (\ref{002}) can be expressed as a hemivariational-inequality.
On the other hand, the notion of a history-dependent operator was introduced in~\cite{SM1} and used in a few recent  papers~\cite{HMS2017,MOS13,MOS18,OGO,SHM,SM2,SMH2018,SP,SX}. For such operator, its current value for a given function at the time instant $t$ depends on the values of the function at the moments from $0$ to $t$. For this reason, we are obviously lead to the so-called history-dependent variational-hemivariational inequalities.
In the present paper, we examine the system  (\ref{001})--(\ref{003}) in the framework of an evolution triple of spaces.

The novelties of the paper are following. 
First, for the first time, we study an abstract  system of two evolution inclusions involving history-dependent operators. 
Until now, 
all aforementioned works have treated only a single quasistatic or evolutionary inclusion, and a single  variational-hemivariational inequality.
As far as we know the system (\ref{001})--(\ref{003}) has not been studied before in the literature.
Here, we fill a gap in the studies and 
provide an existence and uniqueness result 
for a general abstract system  (\ref{001})--(\ref{003}). 

Second, we work under general assumptions on the data which, on the one hand, allow to treat strong coupling between the unknowns, 
and on the other hand guarantee uniqueness 
of solution. Further, 
no compactness hypothesis in the evolution triples in needed as in~\cite{AMMA2}, 
no additional fractional Sobolev spaces are used as in~\cite{KULIG,MIGSZA}, and less restrictive smallness condition is required in comparison with e.g.~\cite{HMS2017,MOSBOOK}.

Third, our results can be applied to several 
particular cases of problem 
(\ref{001})--(\ref{003}) which are of great interest. 
To highlight the generality of the system, 
we study, in Section~\ref{DVHI}, the evolutionary differential variational-hemivariational inequality of the form
%
%
\begin{eqnarray*}
	&&
	\theta'(t) = {\widetilde{A}}\theta(t)
	+{\widetilde{f}}(t,\theta(t),\vartheta w(t))\ \ \,\mbox{\rm a.e.}\ t\in(0,T),
	\\[1mm]
	&&
	\langle w'(t) + {\mathscr N} w(t), v \rangle + {\widetilde{J}}^0(\theta(t), Mw(t); Mv)\ge
	\langle {\widetilde{F}}(t,\theta(t)),v
	\rangle
	\ \ \mbox{for all} \ v \in V, \, \mbox{a.e.} 
	\ t\in (0, T),
	\\[1mm]
	&&w(0)=w_0, \ \ \theta(0)=\theta_0. 
	\end{eqnarray*}
This problem is a special case of the system (\ref{001})--(\ref{003}). 
It 
was studied in~\cite[Problem~1]{JOGO}  
in the context of nonlinear differential hemivariational inequality, and it combines 
a nonlinear evolution equation and a hemivariational inequality of parabolic type.
The existence of solution to this particular form of our system has been proved by the Rothe method in~\cite[Theorem~20]{JOGO} under different hypotheses than in the present paper. We need assumptions on the strong monotonicity of ${\widetilde{A}}$ and the relaxed monotonicity of ${\widetilde{J}}$ 
to guarantee uniqueness, and avoid a compact embedding in an evolution triple of spaces, compactness of the map $\vartheta$, and the compactness of the Nemitsky operator corresponding 
to the map $M$.
Note that there is a vast literature on differential variational inequalities initiated and first systematically discussed 
by Pang-Stewart~\cite{Pang} in finite dimensional spaces. The notion of differential hemivariational inequality was introduced by Liu et al. in~\cite{Liu1} who studied a problem  
consisting of 
a nonlinear evolution equation in a Banach space and a mixed variational quasi hemivariational inequality of elliptic type. 
More details on this topic can be found in~\cite{Gwinner,HanL,Liu1,Stewart} and the references therein.

Fourth, to illustrate the applicability 
of our main result we prove existence and uniqueness of solution to a system of a hyperbolic damped second order variational-hemivariational inequality 
combined with a parabolic hemivariational inequality. This system represents a new mathematical model of a frictional contact problem in thermoviscoelasticity.
In this way, 
we study the thermal effects that are usually neglected in the mathematical analysis even though they have strong influence on the dynamics of the problem. 
It should be mentioned that while the study 
of contact problems is usually achievable 
for many systems without the presence of thermal effects, the nature of the heat equation with a source term depending on the velocity makes the studies very difficult and consequently very interesting and challenging. 
In fact, problem (\ref{001})--(\ref{003}) arises in
a large number of mathematical models which describe dynamic processes of contact
between a deformable body and an obstacle.

%

For a single variational-hemivariational inequality of the form (\ref{001}) with $A$, $J$ and $\varphi$ independent of $\theta$, the main theorem of this paper extends and improves the main result 
of~\cite[Theorem~9]{HMS2017} in several directions.
Further, results on the hemivariational and variational-hemivariational inequalities can be found in monographs~\cite{CLM,Go11.I,Go11.II,HMP,MOSBOOK,
NP,SOFMIG}.

The paper is organized as follows.
Section~\ref{Background} presents the notation and some auxiliary results.
In Section~\ref{inequality} we establish the main result of this paper, Theorem~\ref{SYSTEM}, on the unique solvbility of (\ref{001})--(\ref{003}). 
Section~\ref{DVHI} deals with the differential variational-hemivariational inequality of evolution type.
In Section~\ref{APPL} we examine a dynamic nonsmooth system of thermoviscoelasticity with multivalued and nonmonotone contact and friction conditions for which we prove existence and uniqueness of a weak solution.

\section{Mathematical background and auxiliary results}\label{Background}


Let $X$ be a Banach space, $X^*$ be the topological dual of $X$, and 
$\langle \cdot, \cdot \rangle_{X^* \times X}$ 
be the duality pairing between $X^*$ and $X$. 
Let $\varphi \colon X \to \real \cup \{ +\infty \}$ be a given function. 
The function $\varphi$ is called proper if its effective domain
${\rm dom} \, \varphi = \{ x \in X 
\mid \varphi (x) < +\infty \}\not= \emptyset$.
It is sequentially lower semicontinuous (l.s.c.)
if $x_n \to x$ in $X$ implies
$\varphi (x) \le \liminf \varphi (x_n)$.
If the function $\varphi $ is convex, 
an element $x^* \in X^*$ 
is called a subgradient of $\varphi$ at
$u \in X$, if
\begin{equation*}\label{B.SUBCONV}
\langle x^*, v - u \rangle \le
\varphi(v) - \varphi(u)
\ \ \mbox{\rm for all} \ v \in X.
\end{equation*}
The set of all elements $x^* \in X^*$ which satisfy this inequality is called the convex subdifferential of $\varphi$ at $u$,
and is denoted by $\partial_c \varphi(u)$.
Next, we recall the notion of the
Clarke generalized gradient for a locally
Lipschitz function $\psi \colon X \to \real$, 
see~\cite{Clarke,DMP1,MOSBOOK}.
The generalized gradient of $\psi$
at $u \in X$ is defined by
$$
\partial \psi(u) = \{ u^* \in X^* \mid \psi^{0}(u; v) \ge
\langle u^*, v \rangle \ \ {\rm for\ all\ } \ v \in X \},
$$
where the generalized directional derivative
of $\psi$ at $u \in X$ in the
direction $v \in X$ is given by
$$
\psi^{0}(u; v) =\limsup_{y \to u, \ \lambda \downarrow 0}\frac{\psi(y+ \lambda v) - \psi(y)}{\lambda}.
$$

Let $A \colon X \to X^*$ be a map. 
We say that $A$ is demicontinuous if for all 
$v \in X$, the function 
$u \mapsto \langle A u, v \rangle_{X^* \times X}$
is continuous, i.e., $A$ is continuous as a mapping from $X$ to $X^*$ endowed with $w^*$-topology.
It is hemicontinuous, if for all $u$, $v$, 
$w \in X$, the function
$t \mapsto \langle A(u+tv), w \rangle_{X^* \times X}$ is continuous on $[0,1]$.
$A$ is monotone if $\langle Au - Av, u - v 
\rangle_{X^* \times X} \ge 0$ for all $u$, $v \in X$. 
It is strongly monotone with constant $m > 0$ 
if $\langle Au - Av, u - v \rangle_{X^* \times X} 
\ge m \|u-v\|^2_X$ for all $u$, $v \in X$.
It is known that for a map $A \colon X \to X^*$ 
with $D(A) = X$, the notions of demicontinuity and hemicontinuity coincide, 
see~\cite[Section~1.9, Exercise~I.9]{DMP2}. 
Given a set $D \subset X$, we write  
$\| D \|_X = \sup \{ \| v \|_X \mid v \in D \}$.
%

%
	A triple of spaces 
	$(\mathscr{X}, \mathscr{H}, \mathscr{X^*})$ 
	is said to be an evolution triple, if
	
	\smallskip
	
	\noindent 
	{\rm (a)} \ $\mathscr X$ is a separable reflexive Banach space, $\mathscr X^*$ is its topological dual,
	
	\smallskip
	
	\noindent 
	{\rm (b)} \ $\mathscr H$ is a separable Hilbert space identified with its dual ${\mathscr H}^* \simeq {\mathscr H}$,
	
	\smallskip
	
	\noindent 
	{\rm (c)} \  ${\mathscr X}$ is embedded continuously
	in ${\mathscr H}$,
	denoted ${\mathscr X} \subset {\mathscr H}$, and densely in ${\mathscr H}$.

It is easy to observe that ${\mathscr H}$ is embedded continuously and densely in 
${\mathscr X}^*$, and the duality brackets $\langle \cdot, \cdot 
\rangle_{{\mathscr X}^* \times {\mathscr X}}$
for the pair $({\mathscr X}^*, {\mathscr X})$ and the inner product 
$\langle \cdot, \cdot \rangle_{{\mathscr H}}$ 
on ${\mathscr H}$ 
coincide on ${\mathscr H} \times {\mathscr X}$, 
that is, 
$\langle \cdot, \cdot \rangle_{{\mathscr H} \times {\mathscr X}} = \langle \cdot, \cdot, \rangle_{{\mathscr H}}$.
When no confusion arises, we drop the subscripts.

Next, we recall an existence and uniqueness 
result, see~\cite[Theorem~6]{HMS2017}, 
which is formulated for a subdifferential inclusion  
in the framework of evolution triple of spaces $(\mathscr{X}, \mathscr{H}, \mathscr{X^*})$.
\begin{Problem}\label{Problem1}
	Find $u \in L^2(0, T; {\mathscr X})$ with $u' \in L^2(0,T; {\mathscr X^*})$ such that
	\begin{equation*}
	\begin{cases}
	\displaystyle
	u'(t) + {\mathscr A}(t, u(t)) + \partial \psi (t, u(t)) \ni f(t) \ \ \mbox{\rm a.e.} \ t \in (0,T), \\[1mm]
	u(0) = u_0 .
	\end{cases}
	\end{equation*}
\end{Problem}

Recall that a function 
$u \in L^2(0, T; {\mathscr X})$ with $u'\ \in L^2(0,T; {\mathscr X^*})$ is called a solution 
to Problem~{\rm \ref{Problem1}} 
if there exists
$u^* \in L^2(0,T; {\mathscr X^*})$ such that
$u'(t) + {\mathscr A}(t, u(t)) + u^*(t) = f(t)$, 
$u^*(t) \in \partial \psi (t, u(t))$ 
a.e. $t \in (0,T)$, and $u(0) = u_0$.

\medskip

We need the following hypotheses on the data.

\smallskip
\noindent 
$\underline{H({\mathscr A})}:$ \quad 
$\displaystyle {\mathscr A} \colon (0, T) \times {\mathscr X} \to {\mathscr X}^*$ is such that \smallskip

\lista{
	\item[(a)]  
	${\mathscr A}(\cdot, v)$ is measurable on $(0, T)$ for all $v \in {\mathscr X}$. \smallskip
	\item[(b)]
	${\mathscr A}(t,\cdot)$ is demicontinuous on 
	${\mathscr X}$ for a.e.\ $t \in (0, T)$. \smallskip
	\item[(c)]
	$\| {\mathscr A} (t, v) \|_{{\mathscr X}^*} \le a_0(t) + a_1 \| v \|_{{\mathscr X}}$ for all 
	$v \in {\mathscr X}$,
	a.e.\ $t \in (0,T)$ with $a_0 \in L^2(0,T)$, $a_0 \ge 0$ and $a_1 \ge 0$. \smallskip
	\item[(d)] 
	${\mathscr A}(t, \cdot)$ is strongly monotone for a.e.\ $t \in (0,T)$, i.e., for a constant $m_{{\mathscr A}} > 0$,
	$$\langle {\mathscr A} (t, v_1) 
	- {\mathscr A}(t, v_2), v_1 - v_2 \rangle_{{{\mathscr X}}^* \times 
		{\mathscr X}}
	\ge m_{{\mathscr A}} \| v_1 - v_2 \|_{{\mathscr X}}^2$$
	for all $v_1$, $v_2 \in {\mathscr X}$, a.e.\ $t \in (0, T)$. 
}

\smallskip

\noindent 
$\underline{H(\psi)}:$ \quad 
$\psi \colon (0, T) \times {\mathscr X} \to \real$ is such that 
\smallskip

\lista{
	\item[(a)]  
	$\psi(\cdot, v)$ is measurable on $(0, T)$ for all $v \in {\mathscr X}$. \smallskip
	\item[(b)]
	$\psi(t, \cdot)$ is locally Lipschitz on $V$ for a.e. $t \in (0, T)$. \smallskip
	\item[(c)]
	$\| \partial \psi(t, v) \|_{{\mathscr X}^*} \le  c_0 (t) + c_1 \| v \|_{{\mathscr X}}$ for all $v \in {\mathscr X}$, \smallskip 
	a.e. $t \in (0, T)$ with $c_0 \in L^2(0, T)$, $c_0 \ge 0$, $c_1 \ge 0$. \smallskip
	\item[(d)] 
	$\partial \psi(t, \cdot)$ is relaxed monotone for a.e. $t \in (0,T)$, 
	i.e., for a constant $m_\psi \ge 0$, 
	\[ \langle \partial \psi(t, v_1) - 
	\partial \psi(t, v_2), v_1 - v_2 
	\rangle_{{\mathscr X}^* \times {\mathscr X}} \ge -m_\psi \| v_1 - v_2 \|^2_{\mathscr X} 
	\]
	for all $v_1$, $v_2 \in {\mathscr X}$, 
	a.e.\ $t \in (0, T)$.
}

\smallskip

\noindent 
$\underline{(H_1)}:$ \quad 
$f \in L^2(0, T; {\mathscr X}^*)$, 
$u_0 \in {\mathscr X}$. 

\smallskip

\noindent 
$\underline{(H_2)}:$ \quad $m_{\mathscr A} 
> m_\psi$.

%
%
\begin{Theorem}\label{T1}
	Under hypotheses $H({\mathscr A})$, $H(\psi)$, $(H_1)$ and $(H_2)$, 
	Problem~{\rm \ref{Problem1}} 
	has a unique solution.
\end{Theorem}
The proof of Theorem~\ref{T1} is omitted here. 
It is based on a surjectivity result, 
see e.g. \cite[Theorem 1.3.73]{DMP2}, and  
it was provided in~\cite[Theorem~6]{HMS2017} under more restrictive smallness assumption. However, a careful examination of the proof of a coercivity condition in~\cite{HMS2017} can convinced the reader that the result holds also under a weaker (relaxed) hypothesis $(H_2)$. 

\medskip 

Finally, we state the fixed point result which is a consequence of the Banach contraction 
principle
(see~\cite[Lemma~7]{KULIG} 
or~\cite[Proposition~3.1]{SM2}).
%
\begin{Lemma}\label{CONTR}
	Let ${\mathbb X}$ be a Banach space and $0 < T < \infty$.
	Let $F \colon L^2(0,T; {\mathbb X}) \to L^2(0,T; {\mathbb X})$ be an operator such that
	\begin{equation*}
	\| (F \eta_1)(t) - (F \eta_2)(t)\|^2_{\mathbb X} \le c
	\int_0^t \| \eta_1(s) - \eta_2(s)\|^2_{\mathbb X} \, ds
	\end{equation*}
	for all $\eta_1$, $\eta_2 \in L^2(0,T;{\mathbb X})$, 
	a.e.\ $t \in (0,T)$ with a constant $c > 0$. Then $F$ has a unique fixed point in 
	$L^2(0, T; {\mathbb X})$, 
	i.e., there exists a unique
	$\eta^* \in  L^2(0,T;{\mathbb X})$ such that 
	$F \eta^* = \eta^*$.
\end{Lemma}

\section{System of evolution inclusions}
\label{inequality}

In this section we study the system of first order evolution inclusions with history-dependent operators.

In what follows, 
we need two evolution triples of spaces 
$(V, H, V^*)$ and $(E, X, E^*)$. For each triple, 
we introduce the Bochner spaces on a finite time interval $(0, T)$.
Let 
$$
{\mathbb{W}} = \{ \, w \in L^2(0, T; V) \mid 
w' \in L^2(0,T; V^*) \, \},
$$
where the time derivative $w'$ is understood 
in the distributional sense. 
It is known that
$L^2(0,T;V)^* \simeq L^2(0,T; V^*)$ 
and the space ${\mathbb{W}}$ endowed with the norm
$\| w \|_{{\mathbb{W}}} = \| w \|_{L^2(0, T; V)} + \| w' \|_{L^2(0, T; V^*)}$
is a separable reflexive Banach space, and 
each element of ${\mathbb{W}}$, after a modification on a set of null measure, can be identified with a unique continuous function on $[0, T]$ with values in $H$.
Further, the embedding 
${\mathbb{W}} \subset C(0, T; H)$ is continuous, 
where $C(0, T; H)$ stands for the space of continuous functions on $[0, T]$ 
with values in $H$. 
Analogously, for an evolution triple of spaces $(E, X, E^*)$, we introduce the space
$$
\mathbb{E} = \{ \, \theta \in L^2(0, T; E) \mid 
\theta' \in L^2(0,T; E^*) \, \}.
$$
Moreover, 
let $Q$, $Y$, and $Z$ be Banach spaces.
The system of inclusions under consideration 
reads as follows.
\begin{Problem}\label{SYSTEM} 
	Find $w \in {\mathbb{W}}$ and $\theta \in {\mathbb E}$ such that
	\begin{equation*}
	\begin{cases}
	\displaystyle
	w'(t) + A(t, \theta(t), w(t)) + 
	({R}_1 w)(t) + 
	\partial J (t, \theta(t), ({\mathcal S} w)(t), w(t)) \\[1mm]
	\qquad\qquad\qquad\qquad \quad \ \   
	+ \, \partial_c \varphi (t, \theta(t), 
	({R} w)(t), w(t)) \ni h_1(t) 
	\ \ \mbox{\rm a.e.} \ \ t \in (0,T), \\[1mm]
	\theta'(t) + B(t, w(t), (R_2w)(t), \theta(t))
	+ \partial g(t, w(t), \theta(t)) \ni h_2(t) \ \ \mbox{\rm a.e.} \ \ t \in (0,T), \\[1mm]
	w(0) = w_0, \ \theta(0) = \theta_0 .
	\end{cases}
	\end{equation*}
\end{Problem}

Throughout the paper, the symbols $\partial_c$ and $\partial$ denote the convex subdifferenial and the generalized gradient of Clarke, respectively.
All subdifferentials and generalized gradients 
are always taken with respect to the last variable 
of given functions. 
The operators $R$, $R_1$, $R_2$ and $S$ stand for history-dependent operators.

In the study of Problem~{\rm \ref{SYSTEM}} we need the following hypotheses on the data.

\smallskip
\noindent 
$\underline{H({A})}:$ \quad 
$\displaystyle 
A \colon (0, T) \times X \times V \to V^*$ is such that 
\smallskip

\lista{
	\item[(a)]  
	$A(\cdot, \theta, v)$ is measurable for all $\theta \in X$, 
	$v \in V$. \smallskip
	\item[(b)]
	$A(t,\cdot, v)$ is continuous for all $v \in V$, a.e. $t \in (0, T)$. \smallskip	
	\item[(c)]
	$\| A(t, \theta, v) \|_{V^*} \le a_0(t) + a_1 \| \theta \|_X + a_2 \| v \|_V$ for all $\theta \in X$, $v \in V$, a.e. 
	$t \in (0, T)$ with $a_0 \in L^2(0,T)$, $a_0$, $a_1$, $a_2 \ge 0$. \smallskip
	\item[(d)] 
	$A(t, \theta, \cdot)$ is demicontinuous and there are constants $m_A >0$, 
	${\bar{m}}_A \ge 0$ such that
	$$\langle A (t, \theta_1, v_1) 
	- A(t, \theta_2, v_2), v_1 - v_2 \rangle_{V^* \times V}
	\ge m_A \| v_1 - v_2 \|_{V}^2
	- {\bar{m}}_A \|\theta_1 - \theta_2 \|_X
	\| v_1 - v_2\|_V
	$$
	for all $\theta_1$, $\theta_2 \in X$, 
	$v_1$, $v_2 \in V$, a.e. $t \in (0, T)$. 
}

\smallskip

\smallskip

\noindent 
$\underline{H(J)}:$ \quad 
$J \colon (0, T) \times X \times Z \times V \to \real$ is such that 
\smallskip

\lista{
	\item[(a)]  
	$J(\cdot, \theta, z, v)$ is measurable for all $\theta \in X$, $z \in Z$, $v \in V$. \smallskip
	\item[(b)]
	$J(t, \cdot, \cdot, v)$ is continuous for all $v \in V$, a.e. $t \in (0, T)$. \smallskip
	\item[(c)]
	$J(t, \theta, z, \cdot)$ is locally Lipschitz for all $\theta \in X$, $z \in Z$, a.e. $t \in (0, T)$. 
	\smallskip
	\item[(d)]
	$\| \partial J(t, \theta, z, v) \|_{V^*} \le  c_{0J} (t) + c_{1J} \| \theta\|_X + 
	c_{2J} \| z \|_Z + c_{3J} \| v \|_V$ \smallskip
	for all $\theta \in X$, $z \in Z$, 
	$v \in V$, a.e. $t \in (0, T)$ 
	with $c_{0J} \in L^2(0, T)$, 
	$c_{0J}$, $c_{1J}$, $c_{2J}$, 
	$c_{3J} \ge 0$. \smallskip
	\item[(e)] 
	$J^0(t, \theta_1, z_1, v_1; v_2 - v_1) + J^0(t, \theta_2, z_2, v_2; v_1 - v_2) \\[1mm]
	~~ \qquad \qquad \qquad \qquad 
	\le m_J \, \| v_1 - v_2 \|^2_V + 
	{\bar{m}}_{J}
	\, (\| \theta_1 - \theta_2\|_X 
	+ \| z_1 - z_2 \|_Z ) \| v_1 - v_2 \|_V$ 
	\\[2mm]
	for all $\theta_i\in X$, $z_i \in Z$, 
	$v_i \in V$, 
	$i = 1$, $2$, a.e.\ $t \in (0, T)$ with $m_J$, ${\bar{m}}_J \ge 0$.
}

\smallskip

\noindent 
$\underline{H(\varphi)}:$ \quad 
$\varphi \colon (0, T) \times X \times Y 
\times V \to \real$ is such that 

\smallskip

\lista{
	\item[(a)]  
	$\varphi(\cdot, \theta, y, v)$ is measurable for all $\theta \in X$, $y \in Y$, $v \in V$. \smallskip
	\item[(b)]
	$\varphi(t, \cdot, \cdot, v)$ is continuous for all $v \in V$, a.e. $t \in (0, T)$. \smallskip
	\item[(c)]
	$\varphi(t, \theta, y, \cdot)$ is convex and l.s.c. for all $\theta \in X$, $y \in Y$, a.e. $t \in (0, T)$. \smallskip
	\item[(d)]
	$\| \partial \varphi(t, \theta, y, v) \|_{V^*} \le  c_{0\varphi} (t) +
	c_{1\varphi} \| \theta \|_X 
	+ c_{2\varphi} \| y \|_Y 
	+ c_{3\varphi} \| v \|_V$ 
	for all $\theta$, $y \in Y$, $v \in V$, a.e. $t \in (0, T)$ with $c_{0\varphi} 
	\in L^2(0, T)$, 
	$c_{0\varphi}$, $c_{1\varphi}$, $c_{2\varphi}$, $c_{3\varphi} \ge 0$. \smallskip
	\item[(e)] 
	$
	\varphi (t, \theta_1, y_1, v_2) - 
	\varphi (t, \theta_1, y_1, v_1) + 
	\varphi (t, \theta_2, y_2, v_1) - 
	\varphi (t, \theta_2, y_2, v_2) \\[1mm]
	~~ \qquad \qquad \qquad \qquad
	\le m_{\varphi} \, 
	(\|\theta_1 - \theta_2 \|_X + 
	\| y_1 - y_2 \|_Y ) \| v_1 - v_2 \|_V
	$ \\[2mm]
	for all $\theta_i \in X$, $y_i \in Y$, 
	$v_i \in V$, $i = 1$, $2$, 
	a.e.\ $t \in (0, T)$ with 
	$m_{\varphi} \ge 0$.
}

\smallskip
\noindent 
$\underline{H(B)}:$ \quad 
$\displaystyle 
B \colon (0, T) \times V \times Q \times E \to E^*$ is such that 
\smallskip

\lista{
	\item[(a)]  
	$B(\cdot, v, {\bar{v}}, \theta)$ is measurable for all $v \in V$, ${\bar{v}} \in Q$, $\theta \in E$. \smallskip
	\item[(b)]
	$B(t,\cdot, \cdot, \theta)$ is continuous for all $\theta \in E$, a.e. $t \in (0, T)$. \smallskip	
	\item[(c)]
	$\| B(t, v, {\bar{v}}, \theta) \|_{E^*} 
	\le b_0(t) + b_1 \| v\|_V + 
	b_2 \| {\bar{v}} \|_Q + b_3 \| \theta \|_E$ for all $v \in V$, ${\bar{v}} \in Q$, $\theta \in E$, a.e. $t \in (0, T)$ with $b_0 \in L^2(0,T)$, 
	$b_i \ge 0$, $i=0, \ldots, 3$. \smallskip
	\item[(d)] 
	$B(t, v, {\bar{v}}, \cdot)$ is demicontinuous and there are constants $m_B >0$, 
	${\bar{m}}_B \ge 0$ such that \\[2mm]
	$\quad \langle B (t, v_1, {\bar{v}}_1, \theta_1) 
	- B(t, v_2, {\bar{v}}_2, \theta_2), \theta_1 - \theta_2 \rangle_{E^* \times E}
	\\[1mm]
	~~ \qquad \qquad \qquad
	\ge m_B \| \theta_1 - \theta_2 \|_E^2 
	- {\bar{m}}_B ( \| v_1 - v_2\|_V 
	+ \|{\bar{v}}_1 - {\bar{v}}_2 \|_Q)
	\|\theta_1 - \theta_2 \|_E
	$\\[2mm]
	for all $v_i \in V$, 
	${\bar{v}}_i \in Q$, $\theta_i \in E$, 
	$i=1$, $2$, a.e. $t \in (0, T)$. 
}

\smallskip

\smallskip

\noindent 
$\underline{H(g)}:$ \quad 
$g \colon (0, T) \times V \times E \to \real$ is such that 
\smallskip

\lista{
	\item[(a)]  
	$g(\cdot, v, \theta)$ is measurable for all $v \in V$, $\theta \in E$. \smallskip
	\item[(b)]
	$g(t, \cdot, \theta)$ is continuous for all $\theta \in E$, a.e. $t \in (0, T)$. \smallskip
	\item[(c)]
	$g(t, v, \cdot)$ is locally Lipschitz for all $v \in V$, a.e. $t \in (0, T)$. 
	\smallskip
	\item[(d)]
	$\| \partial g(t, v, \theta) \|_{E^*} \le  c_{0g} (t) + c_{1g} \| v \|_V + 
	c_{2g} \| \theta\|_E$  
	\smallskip
	for all $v \in V$, $\theta \in X$, 
	a.e. $t \in (0, T)$ 
	with $c_{0g} \in L^2(0, T)$, 
	$c_{0g}$, $c_{1g}$, $c_{2g} \ge 0$. \smallskip
	\item[(e)] 
	$g^0(t, v_1, \theta_1; \theta_2 - \theta_1) + g^0(t, v_2, \theta_2; \theta_1 - \theta_2) \\[1mm]
	~~ \qquad \qquad \qquad \qquad 
	\le m_g \, \| \theta_1 - \theta_2 \|^2_E + 
	{\bar{m}}_{g}
	\, \| v_1 - v_2 \|_V 
	\| \theta_1 - \theta_2 \|_E$ 
	\\[2mm]
	for all $v_i \in V$, $\theta_i\in E$, 
	$i = 1$, $2$, a.e.\ $t \in (0, T)$ with $m_g$, ${\bar{m}}_g \ge 0$.
}

\smallskip

\smallskip

\noindent 
$\underline{(H_3)}:$ \quad 
$h_1 \in L^2(0,T; V^*)$, 
$h_2 \in L^2(0,T;E^*)$, 
$w_0 \in V$, $\theta_0 \in E$.

\smallskip

\smallskip

\noindent 
$\underline{(H_4)}:$ \quad 
$m_A > m_J$, $m_B > m_g$.

\smallskip

\smallskip

\noindent 
$\underline{(H_5)}:$ \quad 
${R} \colon L^2(0, T; V) 
\to L^2(0, T; Y)$, 
${R}_1 \colon L^2(0, T; V) 
\to L^2(0, T; V^*)$, 

\qquad
${R}_2 \colon L^2(0, T; V) 
\to L^2(0, T; Q)$,  
${S} \colon L^2(0, T; V) 
\to L^2(0, T; Z)$ are such that

\smallskip

\lista{
	\item[(a)] 
	$\displaystyle 
	\| ({R} v_1)(t) - ({R} v_2)(t) \|_Y \le 
	c_{R} \int_0^t \| v_1(s) - v_2(s) \|_V ds$,
	\\ 
	\item[(b)] 
	$\displaystyle
	\| ({R}_1 v_1)(t) - ({R}_1 v_2)(t) \|_{V^*} \le 
	c_{R_1} \int_0^t \| v_1(s) - v_2(s) \|_V ds$, \\ 
	\item[(c)] 
	$\displaystyle 
	\| ({R}_2 v_1)(t) - ({R}_2 v_2)(t) \|_Q \le 
	c_{R_2} \int_0^t \| v_1(s) - v_2(s) \|_V ds$, \\ 
	\item[(d)] 
	$\displaystyle 
	\| ({S} v_1)(t) - ({S} v_2)(t) \|_Z \le 
	c_{S} \int_0^t \| v_1(s) - v_2(s) \|_V ds$
}
\quad for all $v_1$, $v_2 \in L^2(0, T; V)$, 
a.e.\ $t\in (0, T)$ with $c_R$, $c_{R_1}$, 
$c_{R_2}$, $c_{S} > 0$.


\begin{Remark}\label{RELAXED}
{\rm
	Note that the hypothesis $H(J)${\rm (e)} can be equivalently rewritten as
	\begin{eqnarray}\label{JJJ}
	&&
	\langle v_1^* - v_2^*, v_1-v_1 \rangle_{V^* \times V} \ge 
	-m_J \| v_1 - v_1 \|^2_V \\ [1mm]
	&&
	\qquad \qquad - {\bar{m}}_{J}
	\, (\| \theta_1 - \theta_2\|_X 
	+ \| z_1 - z_2 \|_Z ) \| v_1 - v_2 \|_V
	\nonumber
	\end{eqnarray}
	for all 
	$v_i^* \in \partial J(t, \theta_i, z_i, v_i)$,
	$\theta_i \in X$, $z_i \in Z$, 
	$v_i \in V$, $i = 1$, $2$, a.e.\ $t \in (0, T)$.
	In particular, if $J$ is independent of the variables $\theta$, $z$, the hypothesis $H(J)${\rm (e)} or, equivalently, $(\ref{JJJ})$ reduces to the following relaxed monotonicity condition 
	\begin{equation}\label{FFF}
	\langle v_1^* - v_2^*, v_1 - v_2 \rangle_{V^* \times V} \ge 
	- m_{J} \| v_1 - v_2 \|^2_V
	\end{equation}
	\noindent 
	for all $v_i^* \in \partial J(t, v_i)$, $v_i \in V$, 
	$i = 1$, $2$, a.e. $t \in (0, T)$ with $m_{2J} \ge 0$.
	The latter has been recently used in the literature to study variational-hemivariational inequalities. 
	For examples of nonconvex functions which satisfy the condition $(\ref{FFF})$, we refer to \cite{MOSBOOK}. Further, observe that when 
	$J$ is convex in its last variable, then $(\ref{FFF})$ holds with $m_J = 0$, 
	i.e., the condition $(\ref{FFF})$ simplifies 
	to the monotonicity of convex subdifferential. 
	Evidently, all comments above are also applied to the hypothesis $H(g)${\rm (e)}.
	Furthermore, if $m_J=m_g=0$, then the  hypothesis $(H_4)$ is trivially satisfied.
}
\end{Remark}

The following is the main existence and uniqueness result.

\begin{Theorem}\label{MAIN}
	Under hypotheses $H(A)$, $H(J)$, $H(\varphi)$, $H(B)$, $H(g)$, $(H_3)$--$(H_5)$, 
	Problem~{\rm \ref{SYSTEM}} has a unique solution.
\end{Theorem}
{\bf Proof.} \ 
The proof is carried out in five steps and it is based
on Theorem~\ref{T1} combined with a fixed-point argument of Lemma~\ref{CONTR}.

\smallskip

Let 
$\lambda \in L^2(0,T; X)$, 
$\xi \in L^2(0, T; V^*)$, 
$\eta \in L^2(0, T; Y)$ and 
$\zeta \in L^2(0, T; Z)$ be fixed. 

\noindent 
{\bf Step 1}. 
We define the operator 
${\mathscr A}_\lambda \colon (0,T) 
\times V \to V^*$ and the 
function 
$\psi_{\lambda\xi\eta\zeta} \colon (0, T) 
\times V \to \real$ by 
\begin{eqnarray*}
	&&
	{\mathscr A}_\lambda (t, v) = 
	A(t, \lambda(t), v), \\[1mm]
	&& 
	\psi_{\lambda\xi\eta\zeta} (t, v) = 
	\langle \xi (t), v \rangle_{V^* \times V} 
	+ \varphi (t, \lambda(t), \eta(t), v) 
	+ J(t, \lambda(t), \zeta(t), v) 
\end{eqnarray*}
for all $v \in V$, a.e. $t \in (0,T)$. 
Consider the following auxiliary problem: 
find $w = w_{\lambda\xi\eta\zeta} \in {\mathbb{W}}$ 
such that
\begin{equation}\label{PPP}
\begin{cases}
\displaystyle
w'(t) + {\mathscr A}_\lambda(t, w(t))  
+ \partial \psi_{\lambda\xi\eta\zeta}(t, w(t))
\ni h_1(t) \ \ \mbox{\rm a.e.} \ \ t \in (0,T), \\
w(0) = w_0 .
\end{cases}
\end{equation}

We apply Theorem~\ref{T1} to prove the unique solvability of problem (\ref{PPP}). 
To this end, we  verify the hypotheses of Theorem~\ref{T1}.

From hypotheses $H(A)$(1), (2) we known that $A(\cdot, \cdot, v)$ is a Carath\'eodory function for all $v \in V$. Since $t \mapsto \lambda(t)$ is measurable, from~\cite[Corollary~2.5.24]{DMP1}, we deduce that 
${\mathscr A}_\lambda(\cdot, v)$ is measurable for all $v \in V$, hence $H({\mathscr A})$(a) holds.
The conditions $H({\mathscr A})$(b) and (c) follow 
directly from $H(A)$(d) and (c), respectively. 
By the hypothesis $H(A)$(d) we immediately obtain that ${\mathscr A}(t, \cdot)$ is strongly monotone with constant $m_{\mathscr A} = m_A$ for a.e. $t\in (0, T)$.
Hence, we deduce that $H({\mathscr A})$ is satisfied.

The measurability of $\psi_{\lambda\xi\eta\zeta}(\cdot, v)$ for all $v \in V$ is a consequence of the hypotheses
$H(J)$(a), (b), $H(\varphi)$(a), (b), and the  measurability of functions 
$t \mapsto \lambda(t)$, $t \mapsto \xi(t)$, 
$t \mapsto \eta(t)$ and $t \mapsto \zeta(t)$.
Since $\varphi(t, \theta, y, \cdot)$ is convex and l.s.c.\ for all $\theta \in X$, $y \in Y$, a.e. $t \in (0, T)$, we know by~\cite[Proposition~5.2.10]{DMP1} that
$\varphi(t, \theta, y, \cdot)$ is locally Lipschitz. Due to the condition $H(J)$(c), we deduce that the function 
$\psi_{\lambda\xi\eta\zeta}(t, \cdot)$ is locally Lipschitz on $V$ 
for a.e.\ $t \in (0, T)$, i.e., $H(\psi)$(2) holds.
Recall that $J(t, \theta, z, \cdot)$ and $\varphi(t, \theta, y, \cdot)$ 
are locally Lipschitz for all $\theta \in X$, 
$z \in Z$, $y \in Y$, 
a.e. $t \in (0, T)$, and therefore by~\cite[Proposition~5.6.23]{DMP1}, 
we have
\begin{equation}\label{rule}
\partial \psi_{\lambda\xi\eta\zeta}(t, v) \subset \xi(t) +
\partial_c \varphi (t, \lambda(t), \eta(t), v) + \partial J (t, \lambda(t), \zeta(t), v)  
\end{equation} 
for all $v \in V$ and a.e. $t \in (0,T)$.
Using (\ref{rule}) and the growth conditions $H(J)$(d) and $H(\varphi)$(d), we infer that
$H(\psi)$(c) is satisfied with a nonnegative function $c_{0} \in L^2(0,T)$ and 
$c_1 = \max\{ c_{3\varphi}, c_{3J} \} \ge 0$.

From $H(\varphi)$(c)  and~\cite[Theorem~6.3.19]{DMP1}, it follows that the subdifferential map 
$\partial_c \varphi(t, \theta, y, \cdot)$ is maximal monotone 
for all $\theta \in X$, $y \in Y$, 
a.e. $t \in (0, T)$. Using the monotonicity of this map together with (\ref{JJJ}), 
we obtain
\begin{eqnarray*} 
	&&
	\langle \partial \psi_{\lambda\xi\eta\zeta} (t, v_1) - \partial \psi_{\lambda\xi\eta\zeta} (t, v_2), 
	v_1 - v_2 \rangle_{V^* \times V} \\[1mm]
	&&\qquad 
	= \langle \partial J(t, \lambda(t), \zeta(t), v_1) - \partial J(t, \lambda(t), \zeta(t), v_2), v_1 - v_2 \rangle_{V^* \times V} \\[1mm]
	&&
	\qquad\quad{}
	+\langle \partial_c \varphi(t, \lambda(t), \eta(t), v_1) 
	- \partial_c \varphi(t, \lambda(t), \eta(t), v_2), 
	v_1 - v_2 \rangle_{V^* \times V}\\[1mm]
	&&\qquad 
	\ge - m_{J}\, \| v_1 - v_2 \|^2_V
\end{eqnarray*}
for all $v_1$, $v_2 \in V$, a.e. $t \in (0, T)$. This means that the relaxed monotonicity condition $H(\psi)$(d) is satisfied with 
with $m_{\psi} = m_{J}$. We deduce that the  function 
$\psi_{\lambda\xi\eta\zeta}$ satisfies $H(\psi)$.

Finally, we can readily see that the hypothesis $(H_1)$ follows from $(H_3)$, and the smallness condition $(H_2)$ is a consequence of $(H_4)$.

\smallskip

Having verified the hypotheses of Theorem~\ref{T1}, 
we conclude from this theorem that problem (\ref{PPP}) has a unique solution $w_{\lambda\xi\eta\zeta} \in {\mathbb{W}}$.

\medskip

\noindent 
{\bf Step 2}.
Consider the following problem: find $w= w_{\lambda\xi\eta\zeta} \in {\mathbb{W}}$ 
such that
\begin{equation}\label{*3}
\begin{cases}
\displaystyle
w'(t) + A(t, \lambda(t), w(t)) + \xi(t) +
\partial J (t, \lambda(t), \zeta(t), w(t))
\\[1mm]
\quad \quad \ + \, \partial_c 
\varphi(t, \lambda(t), \eta(t), w(t)) \ni h_1(t)
\ \ \mbox{\rm a.e.} \  t \in (0,T), \\
w(0) = w_0 .
\end{cases}
\end{equation}
It is clear from Step 1 and the inclusion (\ref{rule}) that the problem (\ref{*3}) 
has a solution. 

We will show the uniqueness of solution to (\ref{*3}) and the dependence of this solution 
on $\lambda$, $\xi$, $\eta$ and $\zeta$. 
Let 
$(\lambda_i, \xi_i, \eta_i, \zeta_i) \in 
L^2(0,T; X \times V^* \times Y \times Z)$ 
and $w_i \in {\mathbb{W}}$ for $i=1$, $2$ 
be solutions to the problem (\ref{*3}) corresponding to $(\lambda_i, \xi_i, \eta_i, \zeta_i)$. 
For simplicity of notation, we skip the subscripts 
$\lambda$, $\xi$, $\eta$ and $\zeta$.
We have
\begin{eqnarray*}
	&&
	w_1'(t) + A(t, \lambda_1(t), w_1(t)) + \xi_1(t) +
	\partial J (t, \lambda_1(t), \zeta_1(t), w_1(t))
	\\[1mm]
	&&\quad \quad \, \ + \, \partial_c 
	\varphi(t, \lambda_1(t), \eta_1(t), w_1(t)) \ni h_1(t)
	\ \ \mbox{\rm a.e.} \  t \in (0,T), \\[1mm]
	&&
	w_2'(t) + A(t, \lambda_2(t), w_2(t)) + \xi_2(t) +
	\partial J (t, \lambda_2(t), \zeta_2(t), w_2(t))
	\\[1mm]
	&&\quad \quad \, \ + \, \partial_c 
	\varphi(t, \lambda_2(t), \eta_2(t), w_2(t)) \ni h_2(t)
	\ \ \mbox{\rm a.e.} \  t \in (0,T)
\end{eqnarray*}
with $w_1(0) = w_2(0) = w_0$. We take the difference and multiply in duality by 
$w_1(t) - w_2(t)$ to get
\begin{eqnarray}
&&
\hspace{-0.6cm}
\langle w_1'(t) -w_2'(t), w_1(t) - w_2(t) \rangle_{V^* \times V}\label{***} \\ 
&&
+ \langle 
A(t, \lambda_1(t), w_1(t))  
- A(t, \lambda_2(t), w_2(t)), w_1(t) - w_2(t) \rangle_{V^* \times V} \nonumber\\
&&
+ \langle \xi_1(t) - \xi_2(t),  w_1(t) - w_2(t) \rangle_{V^* \times V} \nonumber\\
&&
+\langle \partial J (t, \lambda_1(t), \zeta_1(t), w_1(t))
- \partial J (t, \lambda_2(t), \zeta_2(t), w_2(t)), w_2(t)), w_1(t) - w_2(t) \rangle_{V^* \times V} \nonumber\\
&& 
+ \langle \partial_c \varphi(t, \lambda_1(t), \eta_1(t), w_1(t)) - \partial_c 
\varphi(t, \lambda_2(t), \eta_2(t), w_2(t)),
w_1(t) - w_2(t) \rangle_{V^* \times V} = 0
\nonumber
\end{eqnarray}
a.e. $t \in (0, T)$.
It is easily seen that the hypothesis $H(\varphi)$(e) is equivalent to the inequality
\begin{eqnarray}\label{*PHI}
&&
\langle \partial \varphi(t, \theta_1, y_1, v_1) - \partial \varphi(t, \theta_2, y_2, v_2),
v_1 - v_2 \rangle_{V^* \times V} \\[1mm]
&&
\qquad\qquad \qquad\quad 
\ge - m_{\varphi} \, 
(\|\theta_1 - \theta_2 \|_X + 
\| y_1 - y_2 \|_Y ) \| v_1 - v_2 \|_V
\nonumber
\end{eqnarray}
for all $\theta_i \in X$, $y_i \in Y$, 
$v_i \in V$, $i = 1$, $2$, 
a.e.\ $t \in (0, T)$.
Moreover, by using the integration by parts formula, 
see~\cite[Proposition~8.4.14]{DMP2},  
and $w_1(0) = w_2(0)$, we have
\begin{equation}\label{byparts}
\int_0^t 
\langle w_1'(s) -w_2'(s), w_1(s) - w_2(s) \rangle_{V^* \times V} \, ds = 
\frac{1}{2} \| w_1(t) - w_2(t) \|^2_H
\end{equation}
for all $t \in [0, T]$.
Combining (\ref{***}) with $H(A)$(d), (\ref{JJJ}), (\ref{*PHI}), and (\ref{byparts}), 
and employing the H\"older inequality, we obtain
\begin{eqnarray*}
	&&
	\hspace{-0.6cm}
	\frac{1}{2} \| w_1(t) - w_2(t) \|^2_H 
	+ (m_A - m_J) \| w_1 - w_2 \|^2_{L^2(0, t; V)}
	\\ [2mm]
	&&
	\le 
	({\bar{m}}_A + m_\varphi + {\bar{m}}_J) 
	\| \lambda_1 - \lambda_2 \|^2_{L^2(0, t; X)}
	\| w_1 - w_2 \|_{L^2(0, t; V)}\\[2mm]
	&& + \, 
	m_\varphi
	\| \eta_1 - \eta_2 \|^2_{L^2(0, t; Y)}
	\| w_1 - w_2 \|_{L^2(0, t; V)} \\[2mm]
	&&
	+ \, 
	{\bar{m}}_J 
	\| \zeta_1 - \zeta_2 \|^2_{L^2(0, t; Z)}
	\| w_1 - w_2 \|_{L^2(0, t; V)}
	+
	\| \xi_1 - \xi_2 \|^2_{L^2(0, t; V^*)}
	\| w_1 - w_2 \|_{L^2(0, t; V)}
\end{eqnarray*}
for all $t \in [0, T]$. 
Taking into account the hypothesis $(H_4)$, we infer that there exists a constant $c>0$ such that
\begin{eqnarray}\label{dep}
&&
\hspace{-0.6cm}
\| w_1 - w_2 \|_{L^2(0, t; V)} \le c \, \Big(
\| \lambda_1 - \lambda_2 \|^2_{L^2(0, t; X)}
+
\| \eta_1 - \eta_2 \|^2_{L^2(0, t; Y)} \\
&&
\qquad\qquad \quad\quad \ \, + \, 
\| \zeta_1 - \zeta_2 \|^2_{L^2(0, t; Z)}
+
\| \xi_1 - \xi_2 \|^2_{L^2(0, t; V^*)} \Big)
\nonumber
\end{eqnarray}
for all $t \in [0, T]$. 
From the inequality (\ref{dep}) it follows that the solution to the problem (\ref{*3}) is unique.  
In conclusion, the problem (\ref{*3}) has a unique solution for which (\ref{dep}) holds.

\medskip

\noindent 
{\bf Step 3}. 
Let $w_{\lambda\xi\eta\zeta} \in {\mathbb{W}}$
be the unique solution to problem (\ref{*3}) obtained in Step~2.
We introduce element
${\bar{w}}_{\lambda\xi\eta\zeta} \in L^2(0,T;Q)$ 
defined by ${\bar{w}}_{\lambda\xi\eta\zeta} = R_2w_{\lambda\xi\eta\zeta}$. 
Consider the following problem:
find $\theta = \theta_{\lambda\xi\eta\zeta} \in {\mathbb E}$ such that
\begin{equation}\label{*5}
\begin{cases}
\displaystyle
\theta'(t) + B(t, w_{\lambda\xi\eta\zeta}(t), 
{\bar{w}}_{\lambda\xi\eta\zeta}(t), \theta(t)) 
+
\partial g (t, w_{\lambda\xi\eta\zeta}(t), \theta(t)) \ni h_2(t)
\ \ \mbox{\rm a.e.} \  t \in (0,T), \\
\theta(0) = \theta_0 .
\end{cases}
\end{equation}
We claim that the problem (\ref{*5}) has the  unique solution such that  
\begin{equation}\label{*6}
\| \theta_1(t) - \theta_2(t)\|^2_X 
+
\int_0^t \| \theta_1(s) - \theta_2(s) \|_{E}^2 \, ds \le 
c \int_0^t \| w_1(s) - w_2(s)\|_V^2 \, ds
\end{equation}
for all $t \in [0, T]$ with a constant $c>0$, where $\theta_1$, $\theta_2 \in {\mathbb E}$ are the unique solutions to (\ref{*5}) corresponding to 
$w_1$, $w_2 \in {\mathbb{W}}$.
To shorten notation, we again skip subscripts $\lambda$, $\xi$, $\eta$ and $\zeta$.

Note that the problem (\ref{*5}) can be formulated as follows: find $\theta \in {\mathbb E}$ such that
\begin{equation*}
\begin{cases}
\displaystyle
\theta'(t) + B_{w,{\bar{w}}}(t,\theta(t)) 
+
\partial \psi_w (t, \theta(t)) \ni h_2(t)
\ \ \mbox{\rm a.e.} \  t \in (0, T), \\
\theta(0) = \theta_0,
\end{cases}
\end{equation*}
where 
$B_{w,{\bar{w}}}\colon (0, T) \times E \to E^*$ 
is defined by 
$$
B_{w,{\bar{w}}} (t, \theta) 
= B(t, w(t),{\bar{w}}(t),\theta)
\ \ \mbox{for} \ \ \theta \in E, 
\ \mbox{a.e.} \ t \in (0,T),
$$ 
and 
$\psi_w \colon (0, T) \times E \to \real$ is given by
$$
\psi_w(t,\theta) = g(t, w(t), \theta)
\ \ \mbox{for} \ \ \theta \in E, 
\ \mbox{a.e.} \ t \in (0, T).
$$
Exploiting the hypotheses $H(B)$ and $H(g)$, we can verify, as in Step 1, that 
$B_{w,{\bar{w}}}$ and $\psi_w$ satisfy the  conditions $H({\mathscr A})$ and $H(\psi)$, 
respectively. Further, the hypotheses $(H_3)$, $(H_4)$ guarantee that $(H_1)$ and $(H_2)$ are satisfied. Applying Theorem~\ref{T1}, we deduce that the problem (\ref{*5}) has a unique solution
$\theta = \theta_{\lambda\xi\eta\zeta} \in {\mathbb E}$. 

It remains to verify (\ref{*6}). 
Let $w_i \in {\mathbb{W}}$ be unique solutions to (\ref{*3}), and 
$\theta_i \in {\mathbb E}$ be the unique solutions to (\ref{*5}) corresponding to 
$w_i \in {\mathbb{W}}$ and 
${\bar{w}}_i = R_2w_i$ for $i=1$, $2$. 
We have
\begin{eqnarray*}
	&&
	\theta_1'(t) + B(t, w_1(t), {\bar{w}}_1(t), \theta_1(t)) 
	+
	\partial g (t, w_1(t), \theta_1(t)) \ni h_2(t), \\[2mm]
	&&
	\theta_2'(t) + B(t, w_2(t), {\bar{w}}_2(t), \theta_2(t)) 
	+
	\partial g (t, w_2(t), \theta_2(t)) \ni h_2(t)
\end{eqnarray*}
for a.e. $t \in (0, T)$, where 
$\theta_1(0) = \theta_2(0) = \theta_0$. 
Similarly, as in Step 2, by using $H(B)$ and $H(g)$, the integration by parts formula, and H\"older's inequality, we obtain
\begin{eqnarray}\label{**99}
&&
\frac{1}{2} \| \theta_1(t) - \theta_2(t) \|_X^2 
+ (m_B - m_g) 
\| \theta_1 - \theta_2\|^2_{L^2(0,t;E)}
\\
&&
\qquad 
\le \Big( ({\bar{m}}_B + {\bar{m}}_g)
\| w_1 - w_2\|_{L^2(0, t; V)}
+ {\bar{m}}_B \| {\bar{w}}_1 - {\bar{w}}_2 \|_{L^2(0,t;Q)} \Big) 
\| \theta_1 - \theta_2\|_{L^2(0,t;E)}
\nonumber
\end{eqnarray}
for all $t \in [0,T]$.
Next, to the right-hand side of (\ref{**99}), 
we apply the Young inequality 
$ab \le \frac{\varepsilon^2}{2} a^2 +
\frac{1}{2\varepsilon^2}b^2$ 
for $a$, $b \in \real$ with 
$\varepsilon^2 = m_B-m_g >0$, 
and estimate it by
\begin{eqnarray*}
	&&
	\frac{m_B-m_g}{2} 
	\| \theta_1 - \theta_2\|^2_{L^2(0,t;E)} \\
	&&\qquad + \, \frac{1}{2(m_B-m_g)}
	\Big( 
	({\bar{m}}_B + {\bar{m}}_g)
	\| w_1 - w_2\|_{L^2(0, t; V)}
	+ {\bar{m}}_B \| {\bar{w}}_1 - {\bar{w}}_2 \|_{L^2(0,t;Q)} \Big)
\end{eqnarray*}
for all $t \in [0,T]$.
On the other hand, by $(H_5)$(c) and the Jensen inequality, we find 
$$
\| {\bar{w}}_1 - {\bar{w}}_2 \|_{L^2(0,t;Q)}
\le c_{R2}\, T \, \| w_1 - w_2 \|_{L^2(0,t;V)}
$$
for all $t \in [0,T]$. 
Using the latter, from (\ref{**99}), we have
\begin{equation*}
\frac{1}{2} \| \theta_1(t) - \theta_2(t) \|_X^2 
+ \frac{m_B - m_g}{2} 
\| \theta_1 - \theta_2\|^2_{L^2(0,t;E)}
\le c
\int_0^t \| w_1(s) -w_2(s) \|^2_V \, ds
\end{equation*}
for all $t \in [0,T]$, where 
$$
c = \frac{1}{m_B-m_g}
\left(
({\bar{m}}_B +{\bar{m}}_g)^2 + {\bar{m}}_B^2 \, c_{R2}^2 \, T^2 \right).
$$
Hence, the estimate (\ref{*6}) follows, which 
completes Step 3. 

\medskip

\noindent 
{\bf Step 4}.
We define the operator 
$F \colon L^2(0,T; X \times V^* \times Y \times Z) \to L^2(0,T; X \times V^* \times Y \times Z)$ by
$$
F(\lambda, \xi, \eta, \zeta) = 
(\theta_{\lambda\xi\eta\zeta}, 
{R}_1 w_{\lambda\xi \eta \zeta}, 
{R} w_{\lambda\xi \eta \zeta}, 
{S} w_{\lambda\xi \eta \zeta}) 
$$
for $(\lambda, \xi, \eta, \zeta) \in 
L^2(0,T; X \times V^* \times Y \times Z)$,
where $w_{\lambda\xi \eta \zeta} 
\in {\mathbb W}$ denotes the unique solution 
to the problem (\ref{*3}), and 
$\theta_{\lambda\xi \eta \zeta} 
\in {\mathbb E}$ is the unique solution 
to the problem (\ref{*5}) corresponding to 
$(\lambda, \xi, \eta, \zeta)$.
We will use Lemma~\ref{CONTR} to show 
that $F$ has a unique fixed point.

Let $(\lambda_i, \xi_i, \eta_i, \zeta_i) \in 
L^2(0,T; X \times V^* \times Y \times Z)$, and
$w_i = w_{\lambda_i\xi_i \eta_i \zeta_i}$, 
$\theta_i = \theta_{\lambda_i\xi_i \eta_i \zeta_i}$, $i=1$, $2$ be the corresponding unique solutions to (\ref{*3}) and (\ref{*5}), respectively. 
By the definition of operator $F$, 
hypothesis $(H_5)$, the estimates (\ref{dep}) and (\ref{*6}), and the Jensen inequality, we have
\begin{eqnarray*}
	&&
	\| F(\lambda_1, \xi_1,\eta_1,\zeta_1)(t) - 
	F(\lambda_2, \xi_2, \eta_2,\zeta_2)(t) 
	\|^2_{X \times V^* \times Y \times Z} \\[2mm]
	&&\qquad = 
	\| \theta_1(t) - \theta_2(t)\|_X^2 +
	\| ({R}_1 w_1)(t) - ({R}_1 w_2)(t) \|_{V^*}^2 \\[2mm]
	&&
	\qquad 
	+ \, \| ({R} w_1)(t) - ({R} w_2)(t) \|_Y^2 
	+ \| ({S} w_1)(t) - ({S} w_2)(t) \|_Z^2 \\[2mm]
	&&\qquad 
	\le 
	c \int_0^t \| w_1(s) - w_2(s) \|^2_V \, ds
	+ \Big( c_{R_1} \int_0^t \| w_1(s) - w_2(s) \|_V \, ds \Big)^2 \\
	&& 
	\qquad 
	+ \, \Big( c_{R} \int_0^t \| w_1(s) - w_2(s) \|_V \, ds \Big)^2 +
	\Big( c_{S} \int_0^t \| w_1(s) - w_2(s) \|_V \, ds \Big)^2 \\ 
	&&\qquad 
	\le c \, \| w_1 - w_2 \|_{L^2(0, t; V)}^2 
	\le c \, \Big(
	\| \lambda_1 - \lambda_2 \|^2_{L^2(0, t; X)}
	+
	\| \xi_1 - \xi_2 \|^2_{L^2(0, t; V^*)}\\[2mm]
	&&\qquad 
	+ \, \| \eta_1 - \eta_2 \|^2_{L^2(0, t; Y)} 
	+ \| \zeta_1 - \zeta_2 \|^2_{L^2(0, t; Z)}
	\Big)
\end{eqnarray*} 
with suitable constants $c>0$. 
This inequality implies
\begin{eqnarray*}
	&&\| F(\lambda_1,\xi_1, \eta_1, \zeta_1)(t) - 
	F(\lambda_1,\xi_2, \eta_2, \zeta_2)(t) 
	\|^2_{X \times V^* \times Y \times Z} 
	\nonumber \\ 
	&&\qquad 
	\le c \, \int_{0}^{t} \, \| (\lambda_1, \xi_1, \eta_1, \zeta_1)(s) - 
	(\lambda_2, \xi_2,\eta_2, \zeta_2)(s) 
	\|^2_{X \times V^* \times Y \times Z} \, ds
	\label{pstaly}
\end{eqnarray*} 
for a.e.\ $t \in (0,T)$. 
We apply Lemma~\ref{CONTR} to deduce that 
there exists a unique fixed point $(\lambda^*,\xi^*,\eta^*,\zeta^*) \in 
L^2(0,T; X \times V^* \times Y \times Z)$ of $F$, i.e.,
$$
F(\lambda^*, \xi^*, \eta^*, \zeta^*) 
= (\lambda^*, \xi^*, \eta^*, \zeta^*).
$$

\smallskip

\noindent 
{\bf Step 5.} 
Let $(\lambda^*, \xi^*,\eta^*, \zeta^*) 
\in L^2(0,T; X \times V^* \times Y \times Z)$ 
be the unique fixed point of the operator $F$. 
Let $w^* = w_{\lambda^* \xi^* \eta^* \zeta^*} 
\in {\mathbb{W}}$ be the unique solution 
to the problem (\ref{*3}) 
corresponding to $(\lambda^*, \xi^*,\eta^*, \zeta^*)$, see Step~2.

Next, define ${\bar{w}}^* \in L^2(0,T; Q)$ by 
${\bar{w}}^* = R_2 w^*$. 
Further, let $\theta^* 
= \theta_{\lambda^* \xi^* \eta^* \zeta^*} 
\in {\mathbb E}$ be the unique solution 
to the problem (\ref{*5}) 
which corresponds to $w^*$ and ${\bar{w}}^*$, 
see Step~3. 
By the definition of operator $F$, we have
$$
\lambda^* = \theta^*, \ \ 
\xi^* = R_1 w^*, \ \ 
\eta^* = Rw^*
\ \ \mbox{and} \ \ 
\zeta^* = S w^*.  
$$
Using these relations in (\ref{*3}) and (\ref{*5}), 
we conclude that $(w^*,\theta^*) \in {\mathbb{W}} \times {\mathbb{E}}$ is the unique solution to Problem~\ref{SYSTEM}. This completes the proof.
\hfill$\Box$

\begin{Remark}\label{Remark22}
{\rm
The following conditions are useful to check the hypotheses $H(A)${\rm (d)}  
and $H(B)${\rm (d)}.
If $A(t,\theta, \cdot)$
is strongly monotone with $m_A>0$ 
for all $\theta \in X$, 
a.e. $t \in (0,T)$,  
and $A(t,\cdot, v)$ is 
Lipschitz with $L_A>0$ 
for all $v \in V$, a.e. $t \in (0, T)$, 
i.e., 
\begin{eqnarray*}
&&
\langle A (t, \theta, v_1) 
- A(t, \theta, v_2), v_1 - v_2 
\rangle_{V^* \times V}
\ge m_A \| v_1 - v_2 \|_{V}^2, \\[2mm]
&&
\| A (t, \theta_1, v) 
- A(t, \theta_2, v) \|_{V^*} \le L_A\, 
\| \theta_1 - \theta_2\|_X,
\end{eqnarray*}
then
\begin{equation*}
\langle A (t, \theta_1, v_1) 
- A(t, \theta_2, v_2), v_1 - v_2 \rangle_{V^* \times V}
\ge m_A \, \| v_1 - v_2 \|_{V}^2
- L_A \, \|\theta_1 - \theta_2 \|_X
\| v_1 - v_2\|_V
\end{equation*}
for all $\theta_1$, $\theta_2 \in X$, 
$v_1$, $v_2 \in V$, a.e. $t \in (0, T)$. 
}
\end{Remark}

\section{Differential variational-hemivariational inequality}\label{DVHI}

In this section we provide the first application of Theorem~\ref{MAIN} to study the evolutionary differential variational-hemivariational inequality.
%
\begin{Problem}\label{DVHVI} 
	Find $u \in {\mathbb{W}}$ and $\theta \in {\mathbb E}$ such that
	\begin{equation*}
	\begin{cases}
	\displaystyle
	\langle u'(t) + {\bar{A}}(t, u(t)), v - u(t) \rangle_{V^* \times V} +
	G^0(t,\theta(t), Mu(t); Mv - Mu(t))) +
	\\[1mm]
	\ \   
	+ \, \varphi (t, \theta(t), v) - \varphi(t,\theta(t),u(t)) \ge 
	\langle F(t, \theta(t)), v - u(t) \rangle_{V^* \times V} \ \, \mbox{\rm for all} \ v \in V, 
	\ \mbox{\rm a.e.} \ t \in (0,T), \\[1mm]
	\theta'(t) + {\bar{B}}(t, \theta(t)) =
	f (t, \theta(t), \vartheta u(t))) 
	\ \ \mbox{\rm a.e.} \ t \in (0,T), \\[1mm]
	u(0) = u_0, \ \theta(0) = \theta_0 .
	\end{cases}
	\end{equation*}
\end{Problem}

This system represents the differential variational-hemvariational inequality and 
consists of the variational-hemivariational 
inequality of parabolic type driven 
by the nonlinear abstract evolution equation. 
To the best of our knowledge, a simpler version of Problem~\ref{DVHVI} was studied for the first time in~\cite{JOGO}, only in a special case with  $\varphi=0$. 
There, the solution to the evolution equation was considered in a mild sense with the time-independent operator ${\bar{B}}$ being the infinitesimal generator of a $C_0$-semigroup. In~\cite{JOGO}, we have established existence of solution by the Rothe method. 
Below, we study Problem~\ref{DVHVI} 
which is more general than that one in~\cite{JOGO},
by using an alternative framework, 
under the different hypotheses than those in~\cite{JOGO}, and we prove the existence and uniqueness of solution. We will use the results of Section~\ref{inequality} and, for simplicity, 
we consider a problem without history-dependent operators.

Let $(V,H,V^*)$ and $(E, X, E^*)$ be two evolution triples of spaces, and $X_0$ and $Y_0$ be reflexive and separable Banach spaces. We need the following hypotheses.

\smallskip
\noindent 
$\underline{H({\bar{A}})}:$ \quad 
$\displaystyle 
{\bar A} \colon (0, T) \times V \to V^*$ is such that 
\smallskip

\lista{
	\item[(a)]  
	${\bar A}(\cdot, v)$ is measurable for all  
	$v \in V$. \smallskip
	\item[(b)]
	$\| {\bar A}(t, v) \|_{V^*} \le {\bar a}_0(t) + {\bar a}_1 \| v \|_V$ for $v \in V$, a.e. 
	$t \in (0, T)$ with ${\bar a}_0 \in L^2(0,T)$, ${\bar a}_0$, ${\bar a}_1 \ge 0$. \smallskip
	\item[(c)] 
	${\bar A}(t, \cdot)$ is demicontinuous 
	and strongly monotone with constant 
	$m_{\bar A} >0$, a.e. $t \in (0, T)$.  
}

\smallskip

\smallskip

\noindent 
$\underline{H(M, \vartheta)}:$ \quad 
$M \colon V \to X_0$ and 
$\vartheta \colon H \to Y_0$ 
are linear and bounded.
 
\smallskip

\smallskip

\noindent 
$\underline{H(G)}:$ \quad 
$G \colon (0, T) \times X \times X_0 \to \real$ is such that 

\smallskip

\lista{
	\item[(a)]  
	$G(\cdot, \theta, v)$ is measurable for all $\theta \in X$, $v \in X_0$. \smallskip
	\item[(b)]
	$G(t, \cdot, v)$ is continuous for all 
	$v \in X_0$, a.e. $t \in (0, T)$. \smallskip
	\item[(c)]
	$G(t, \theta, \cdot)$ is locally Lipschitz for all $\theta \in X$, a.e. $t \in (0, T)$. 
	\smallskip
	\item[(d)]
	$\| \partial G(t, \theta, v) \|_{X_0^*} \le  
	c_{0G} (t) + c_{1G} \| \theta\|_X 
	+ c_{2G} \| v \|_{X_0}$ \smallskip
	for all $\theta \in X$, 
	$v \in X_0$, a.e. $t \in (0, T)$ 
	with $c_{0G} \in L^2(0, T)$, 
	$c_{0G}$, $c_{1G}$, $c_{2G} \ge 0$. \smallskip
	\item[(e)] 
	$G^0(t, \theta_1, v_1; v_2 - v_1) + 
	G^0(t, \theta_2, v_2; v_1 - v_2) \\[1mm]
	~~ \qquad \qquad \qquad \qquad 
	\le m_G \, \| v_1 - v_2 \|^2_{X_0} + 
	{\bar{m}}_{G}
	\, \| \theta_1 - \theta_2\|_X 
	\| v_1 - v_2 \|_{X_0}$ 
	\\[2mm]
	for all $\theta_i\in X$, $v_i \in X_0$, 
	$i = 1$, $2$, a.e.\ $t \in (0, T)$ with $m_G$, ${\bar{m}}_G \ge 0$.
}

\smallskip

\smallskip

\noindent 
$\underline{H(\varphi)_1}:$ \quad 
$\varphi \colon (0, T) \times X \times V \to \real$ 
satisfies $H(\varphi)$ without the $y$-variable. 

\smallskip

\smallskip

\noindent 
$\underline{H(F)}:$ \quad 
$F \colon (0, T) \times X \to V^*$ is such that 

\smallskip

\lista{
	\item[(a)]  
	$F(\cdot, \theta)$ is measurable for all $\theta \in X$. \smallskip
	\item[(b)]
	$F(t, \cdot)$ is Lipschitz continuous with constant $L_F > 0$, a.e. $t \in (0, T)$. \smallskip
	\item[(c)]
	$\| F(t, \theta) \|_{V^*} \le  
	c_{0F} (t) + c_{1F} \| \theta\|_X$ 
	\smallskip
	for all $\theta \in X$, a.e. $t \in (0, T)$ 
	with $c_{0F} \in L^2(0, T)$, 
	$c_{0F}$, $c_{1F} \ge 0$.
}

\smallskip
\noindent 
$\underline{H({\bar{B}})}:$ \quad 
$\displaystyle 
{\bar B} \colon (0, T) \times E \to E^*$ is such that 
\smallskip

\lista{
	\item[(a)]  
	${\bar B}(\cdot, \theta)$ is measurable for all  
	$\theta \in E$. \smallskip
	\item[(b)]
	$\| {\bar B}(t, \theta) \|_{E^*} \le 
	{\bar b}_0(t) + {\bar b}_1 \| \theta \|_E$ 
	for $\theta \in E$, a.e. 
	$t \in (0, T)$ with ${\bar b}_0 \in L^2(0,T)$, ${\bar b}_0$, ${\bar b}_1 \ge 0$. \smallskip
	\item[(c)] 
	${\bar B}(t, \cdot)$ is demicontinuous 
	and strongly monotone with constant 
	$m_{\bar B} >0$, a.e. $t \in (0, T)$.  
}

\smallskip

\smallskip

\noindent 
$\underline{H(f)}:$ \quad 
$f \colon (0, T) \times E \times Y_0 \to E^*$ 
is such that 

\smallskip

\lista{
	\item[(a)]  
	$f(\cdot, \theta, v)$ is measurable for all $\theta \in E$, $v \in Y_0$. \smallskip
	\item[(b)]
	$f(t, \theta, \cdot)$ is Lipschitz 
	continuous with $L_f > 0$ 
	for all $\theta \in E$, a.e. $t \in (0, T)$. \smallskip
	\item[(c)]
	$\| f(t, \theta, v) \|_{E^*} \le  
	c_{0f} (t) + c_{1f} \| \theta\|_E 
	+ c_{2f} \| v \|_{Y_0}$ \smallskip
	for all $\theta \in E$, 
	$v \in Y_0$, a.e. $t \in (0, T)$ 
	with $c_{0f} \in L^2(0, T)$, 
	$c_{0f}$, $c_{1f}$, $c_{2f} \ge 0$. \smallskip
	\item[(e)] 
	$f(t,\cdot, v)$ is demicontinuous and one-sided Lipschitz continuous, i.e.,
	$$
	\langle f(t,\theta_1, v) - f(t,\theta_2,v),
	\theta_1 - \theta_2 \rangle_{E^* \times E}
	\le L_{os} \| \theta_1 - \theta_2 \|^2_E
	$$
	for all $\theta_1$, 
	$\theta_2 \in E$, $v \in Y_0$  
	with $L_{os} > 0$, a.e. $t \in (0, T)$.
}

\smallskip

\smallskip

\noindent 
$\underline{(H_6)}:$ \quad 
$u_0 \in V$, \ $\theta_0 \in E$.

\smallskip

\smallskip

\noindent 
$\underline{(H_7)}:$ \quad 
$m_{\bar A} > m_G \| M \|^2$, \ 
$m_{\cal B} > L_{os}$.

\smallskip

\begin{Theorem}\label{T2}
Under hypotheses $H({\bar A})$, $H(M,\vartheta)$, $H(G)$, $H(\varphi)_1$, $H(F)$, $H({\bar B})$, $H(f)$, $(H_6)$, and $(H_7)$, 
Problem~{\rm \ref{DVHVI}} has the unique solution
$(u, \theta) \in {\mathbb{W}} \times {\mathbb E}$.
\end{Theorem}
{\bf Proof.} \ 
First, using the estimates obtained in Theorem~\ref{MAIN}, we show the uniqueness of solution.
Let $(u_i, \theta_i) \in {\mathbb{W}} 
\times {\mathbb E}$ be solutions to 
Problem~{\rm \ref{DVHVI}, $i=1$, $2$. 
We have
\begin{eqnarray}
&&\hspace{-1.0cm}
\langle u_i'(t) + {\bar{A}}(t, u_i(t)), v - u_i(t) \rangle_{V^* \times V} +
G^0(t,\theta_i(t), Mu_i(t); Mv - Mu_i(t))) +
\nonumber
\\[1mm]
&&\hspace{-1.0cm}
\ \   
+ \, \varphi (t, \theta_i(t), v) - \varphi(t,\theta_i(t),u_i(t)) \ge 
\langle F(t, \theta_i(t)), 
v - u_i(t) \rangle_{V^* \times V} \ \mbox{\rm for all} \ v \in V, 
\ \mbox{\rm a.e.} \ t , \label{QQ1}
\\[1mm]
&&\hspace{-1.0cm}
\theta_i'(t) + {\bar{B}}(t, \theta_i(t)) =
f (t, \theta_i(t), \vartheta u_i(t))) 
\ \ \mbox{\rm a.e.} \ t \in (0,T), \label{QQ2} 
\\[1mm]
&&\hspace{-1.0cm}
u(0) = u_0, \ \theta(0) = \theta_0 . \nonumber 
\end{eqnarray}
\noindent 
From (\ref{QQ1}), we obtain
\begin{eqnarray*}
&&\hspace{-0.7cm}
\langle u_1'(t) -u_2'(t), u_1(t) - u_2(t) \rangle_{V^* \times V} 
+ 
\langle {\bar{A}}(t, u_1(t)) - {\bar{A}}(t, u_2(t)),
u_1(t) - u_2(t) \rangle_{V^* \times V} \\[1mm]
&&
+\, 
G^0(t,\theta_1(t), Mu_1(t); M u_2(t)- M u_1(t))) +
G^0(t,\theta_2(t), Mu_2(t); Mu_1(t) - Mu_2(t)))
\\
&&
+ \, 
\varphi (t, \theta_1(t), u_2(t)) - \varphi(t,\theta_1(t), u_1(t))
+
\varphi (t, \theta_2(t), u_1(t)) - \varphi(t,\theta_2(t), u_2(t)) \\[1mm]
&&
\le 
\langle F(t, \theta_2(t)) - F(t, \theta_1(t)), 
u_2(t) - u_1(t) \rangle_{V^* \times V} 
\end{eqnarray*}
for a.e. $t \in (0, T)$.
Integrating on $(0, t)$, by the hypotheses, 
we have 
\begin{eqnarray*}
&&\hspace{-0.7cm}
\frac{1}{2} \| u_1(t) - u_2(t) \|_H^2
+ 
\left( m_{\bar{A}} - m_G \| M \|^2	\right)
\int_0^t \| u_1(s) - u_2(s) \|_V^2 \, ds \\
&&\qquad 
\le
\left(
{\bar{m}}_G \| M \| + L_F + m_\varphi \right)
\int_0^t \| \theta_1(s)- \theta_2(s)\|_X
\| u_1 (s) - u_2(s) \|_V \, ds
\end{eqnarray*}
for all $t \in [0,T]$, which implies 
\begin{equation}\label{*6}
\| u_1 - u_2 \|_{L^2(0,t;V)} \le 
c \, \| \theta_1 - \theta_2 \|_{L^2(0,t;X)}
\ \ \mbox{\rm for all} \ t \in [0, T]
\end{equation}
with a constsnt $c > 0$.
From (\ref{QQ2}), we obtain
\begin{eqnarray*}
	&&\hspace{-0.7cm}
	\langle \theta_1'(t) - \theta_2'(t), \theta_1(t) - \theta_2(t) \rangle_{E^* \times E} 
	+ 
	\langle {\bar{B}}(t, \theta_1(t)) - {\bar{B}}(t, \theta_2(t)),
	\theta_1(t) - \theta_2(t) \rangle_{E^* \times E} \\[1mm]
	&&\qquad \qquad\qquad \ \ 
	=\langle
	f(t,\theta_1(t), \vartheta u_1(t)) -
	f(t,\theta_2(t), \vartheta u_2(t)),
	\theta_1(t) - \theta_2(t) \rangle_{E^* \times E}
\end{eqnarray*}
for a.e. $t \in (0, T)$.
Exploiting the hypotheses, we get 
\begin{eqnarray*}
	&&\hspace{-0.7cm}
	\frac{1}{2} \| \theta_1(t) - \theta_2(t) \|_X^2
	+ 
	\left( m_{\bar{B}} - L_{os}\right)
	\int_0^t \| \theta_1(s) - \theta_2(s) \|_E^2 \, ds \\
	&&\qquad 
	\le
	L_f \, \| \vartheta \|
	\int_0^t \| u_1 (s) - u_2(s) \|_V 
	\| \theta_1(s)- \theta_2(s)\|_E  \, ds
\end{eqnarray*}
for all $t \in [0,T]$. Applying the Young inequality (see the proof of Theorem~\ref{MAIN}) 
with $\varepsilon^2 = m_{\bar{B}} - L_{os} > 0$ to the right-hand side, we deduce 
\begin{equation*}
\| \theta_1(t) - \theta_2(t) \|_X^2
+ \| \theta_1 - \theta_2 \|^2_{L^2(0,t;E)} 
\le c \, \| u_1 - u_2 \|^2_{L^2(0,t;V)}
\end{equation*}
for all $t \in [0, T]$ with $c > 0$.
Combining the latter with (\ref{*6}) gives
\begin{equation*}
\| \theta_1(t) - \theta_2(t) \|_X^2
\le c \, \int_0^t \| \theta_1(s) - \theta_2(s) \|^2_X 
\, ds
\ \ \mbox{\rm for all} \ t \in [0, T], 
\end{equation*}
which by the Gronwall inequality implies 
$\theta_1 = \theta_2$. Finally, from (\ref{*6}),  
we have $u_1 = u_2$. This completes the proof of uniqueness to Problem~\ref{DVHVI}.

For the existence part, we will use Theorem~\ref{MAIN}. 
Let $A \colon (0, T) \times X \times V \to V^*$ 
and 
$J \colon (0, T) \times X \times Z \times V \to \real$ be defined by
\begin{eqnarray*}
&&
A(t, \theta, v) = {\bar A}(t, v)
\ \ \mbox{\rm for} \ \ \theta \in X, \, v \in V, \ 
\mbox{a.e.} \ t \in (0,T), \\ [1mm]
&&
J(t, \theta, z, v) = G(t, \theta, Mv) 
- \langle F(t, \theta), v \rangle_{V^*\times V}
\ \ \mbox{\rm for} \ \ \theta \in X, \, v \in V,
\ \mbox{a.e.} \ t \in (0,T), 
\end{eqnarray*}
respectively.
It is easy to check that operator $A$ satisfies $H(A)$ with $m_A = m_{\bar A}$ and ${\bar{m}}_A = 0$.
Since the second term in the definition of $J$ is strictly differentiable, by using~\cite[Proposition~3.37(ii)]{MOSBOOK}, we have 
\begin{equation}\label{*333}
\partial J(t,\theta, z, v) = 
\partial (G(t,\theta, Mv)) - F(t,\theta) 
\subset M^* \partial G(t,\theta,Mv) - F(t, \theta).
\end{equation}
All conditions in $H(J)$ follow from the hypotheses. 
In particular, $H(J)$ holds with 
$m_J = m_G \| M \|^2$ and ${\bar{m}}_J = 
{\bar{m}}_G \| M \| + L_F$.
Next, we define 
$B \colon (0, T) \times V \times Q \times E \to E^*$
by 
$$
B(t,v,{\bar{v}}, \theta)= {\bar{B}}(t, \theta) - 
f(t, \theta, \vartheta v).
$$
for $v \in V$, \, $\theta \in E$, \,  
a.e. $t \in (0,T)$.
It can be verified that $H(B)$ holds with 
$m_B = m_{\bar B}-L_{os} >0$ and ${\bar{m}}_B = 
L_f \| \vartheta \|$.
Next, we choose $g = 0$, $h_1=0$ and $h_2 =0$.
The condition $(H_6)$ entails $(H_3)$, 
$H(\varphi)_1$ implies $H(\varphi)$, 
and $(H_4)$ is a consequence of $(H_7)$. 
The hypothesis $(H_5)$ holds trivially.

Applying Theorem~\ref{MAIN}, we deduce that there exists a unique solution to the following system:
\begin{Problem}\label{DVHVI2}
Find $(u, \theta) \in {\mathbb{W}} \times {\mathbb E}$ such that
\begin{equation*}
	\begin{cases}
	\displaystyle
	u'(t) + {\bar{A}}(t, u(t)) +
	\partial (G(t, \theta(t), Mu(t))) 
	+ \, \partial_c \varphi (t, \theta(t), u(t)) 
	\ni F(t,\theta(t)) 
	\ \ \mbox{\rm a.e.} \ \, t \in (0,T), \\[1mm]
	\theta'(t) + {\bar{B}}(t, \theta(t))
	= f(t, \theta(t), \vartheta u(t))
	\ \ \mbox{\rm a.e.} \ \ t \in (0,T), \\[1mm]
	u(0) = w_0, \ \theta(0) = \theta_0 .
	\end{cases}
\end{equation*}
\end{Problem}
\noindent 
It follows from inclusion (\ref{*333}) and the definition of generalized gradient that 
$(u, \theta) \in {\mathbb{W}} \times {\mathbb E}$ is  also a solution to Problem~\ref{DVHVI}. 
This completes the proof of  existence of solution 
to Problem~\ref{DVHVI}.
Since the solution to Problem~\ref{DVHVI} is unique,  it is clear that Problems~\ref{DVHVI} and~\ref{DVHVI2}
are equivalent, and they have the same unique solution. This completes the proof.
\hfill$\Box$

\section{A frictional contact problem in thermoviscoelasticity}\label{APPL}

In this section we provide the second application of 
our abstract results of Section~\ref{inequality} 
to examine a mathematical model of a dynamic thermoviscoelastic problem which describes 
nonsmooth 
frictional contact between a body and a foundation. 
The model consists of the coupled system of a damped hemivariational inequality of hyperbolic type for the displacement field and a parabolic hemivariational inequality for the temperature. 
We will show by applying Theorem~\ref{MAIN} 
that the contact problem has a unique weak solution.

Recently dynamic frictional contact problems with or without thermal effects for viscoelastic bodies 
have been investigated in a large number of papers, see e.g.~\cite{ADLY,HS,KULIG,KS2001,
MIGORSKI7,MOSBOOK} 
and the references therein.
In spite of importance of the subject in applications, to the best
of the authors’ knowledge, the existence of solutions to the system of hemivariational
inequalities in dynamic thermoviscoelasticity has studied in a few papers~\cite{DM1,DM2,DM3,MIGSZA}.
In these papers, there is a coupling only between the displacement (and velocity) and
the temperature in the constitutive law which is assumed to be linear. In contrast to
the aforementioned papers, now we deal with a nonlinear constitutive relation and assume that the coupling appears not only in the constitutive
law but also in the heat flux boundary condition on the contact surface. 
For more detailed justification,  description, and discussion on mechanical interpretation, we refer to~\cite{HS,MOSBOOK,SST}.
We need an additional notation to formulate the problem.

Let $\Omega \subset \real^d$, $d = 2$, $3$, 
be an open bounded set 
with a Lipschitz boundary 
$\Gamma$. 
The latter consists of parts 
sets $\overline{\Gamma}_D$, $\overline{\Gamma}_N$ and
$\overline{\Gamma}_C$, with mutually disjoint relatively open sets
$\Gamma_D$, $\Gamma_N$ and $\Gamma_C$, such that 
${\rm m} \, (\Gamma_D) >0$.  
A viscoelastic body in its reference configuration occupies volume $\Omega$, 
it is assumed to be stress free 
and at a constant temperature, 
conveniently set as zero. 
We are interested in a mathematical model that describes the evolution 
of the mechanical state of the body and its temperature during 
the finite time interval $[0, T]$. 
The unknows in the problem are 
$\bu \colon {\Omega}\times[0,T] \to \mathbb{R}^d$, $\bu=(u_i)$
$\bsigma \colon {\Omega}\times[0,T] \to \mathbb{S}^d$, $\bsigma = (\sigma_{ij})$ and 
$\theta \colon \Omega \times [0, T] \to \real$ 
which represent the displacement field, the stress field and the temperature, respectively.
Here $\mathbb{S}^{d}$ is the space of second order symmetric tensors endowed with the canonical inner product
$\bsigma\cdot\btau =\sigma_{ij}\tau_{ij}$,
for $\bsigma=(\sigma_{ij})$, $\btau=(\tau_{ij})
\in\mathbb{S}^{d}$.
We often suppress the explicit dependence 
of the quantities on the spatial variable $\bx$, or both $\bx$ and~$t$. 

We suppose that the body is clamped on $\Gamma_D$, the volume forces of density 
$\fb_0$ act in $\Omega$ and 
the surface tractions of density $\fb_N$ 
are applied on $\Gamma_N$. 
The body is subjected to a heat source term per unit volume $g$ and it comes in contact with an obstacle, the so-called foundation, over the contact surface $\Gamma_C$. 
We use the notation 
$Q_T = \Omega \times (0, T)$,   
$\Sigma_D = \Gamma_D \times (0, T)$,
$\Sigma_N = \Gamma_N \times (0, T)$ 
and 
$\Sigma_C = \Gamma_C \times (0, T)$. 
Moreover, 
we denote by ${\bnu}=(\nu_i)$ 
the outward unit normal on the boundary, and 
by $\bvarepsilon(\bu)=(\varepsilon_{ij}(\bu))$ the linearized strain tensor whose components  are given by
\begin{equation*}
\label{Ddef} \varepsilon_{ij}(\bu)= \frac 12\, (u_{i,j} + u_{j,i})
\end{equation*}
where $u_{i,j}=\partial u_i/\partial x_j$. For a vector field, we
use the notation $v_\nu$ and $\bv_\tau$ for the normal and
tangential components of $\bv$ on $\Gamma$ given by
$v_\nu=\bv\cdot\bnu$ and $\bv_\tau=\bv-v_\nu\bnu$. The normal and
tangential components of the stress field $\bsigma$ on the
boundary are defined by $\sigma_\nu=(\bsigma\bnu)\cdot\bnu$ and
$\bsigma_\tau=\bsigma\bnu-\sigma_\nu\bnu$, respectively. 
The classical formulation of the mechanical problem of frictional contact for the thermoviscoelastic body is the following. 
\begin{Problem}\label{p222} 
Find a displacement field 
$\bu \colon Q_T \to\mathbb{R}^d$, 
a stress field
$\bsigma \colon Q_T \to \mathbb{S}^d$ and a temperature 
$\theta \colon Q_T \to \real$
such that for all $t\in (0,T)$,  
\begin{eqnarray} \label{EQ1} 
\bu'' (t) - {\rm Div} \, \bsigma (t) = \fb_0(t) & \mbox{\rm in} & Q_T 
\\
\bsigma(t) = 
{\mathscr A}(t,\theta(t), \bvarepsilon(\bu'(t))) +
{\mathscr B}(t,\bvarepsilon(\bu(t))) + \qquad \qquad 
\nonumber 
\\
+ \int_0^t\,{\mathscr C}(t-s) \bvarepsilon(\bu'(s)) \, ds + 
{\mathscr C}_e (t, \theta (t)) & \mbox{\rm in} & Q_T \label{EQ2}
\\ 
\label{EQ3}  
\theta' (t) - {\rm div} \, {\mathscr K}(t, \nabla \theta (t)) = {\mathscr N}(t, \bu'(t)) + h_0(t) & \mbox{\rm in} & Q_T 
\\
\label{EQ4}  
\bu (t) = 0 & \mbox{\rm on} & \Sigma_D 
\\ 
\label{EQ5}  
\bsigma (t) \bnu = \fb_N (t) &  \mbox{\rm on} & \Sigma_N
\\ 
\label{EQ6}  
-\sigma_\nu(t) \in 
k(t, {\mathscr R} \theta(t), u_\nu(t))\,
\partial j_\nu({u}_\nu'(t))
& \mbox{\rm \rm on} & \Sigma_C   
\\ 
\label{EQ7} 
\|\bsigma_\tau(t)\| \le F_b\Big(t, {\mathscr R} \theta(t), \int_0^t\|\bu_\tau(s)\|\,ds\Big)\ &{\rm on}\ &\Sigma_C \label{H66}
\\
-\bsigma_\tau=F_b\Big(t, {\mathscr R} \theta(t),
\int_0^t\|\bu_\tau(s)\|\,ds\Big)
\frac{{\bu}_\tau'(t)}{\|{\bu}_\tau'(t)\|}\ \ {\rm  if}\
{\bu}_\tau'(t)\ne\bzero \ &{\rm on}\ &\Sigma_C\label{EQ8}
\\
- \frac{\partial \theta}{\partial \nu_{\mathscr K}} 
\in \partial j(t, \theta(t)) - h_\tau(\|\bu_\tau'\|)
\ &{\rm on}\ &\Sigma_C\label{flux}
\\
\label{EQ9}   
\theta (t) = 0 & {\rm on} & \Sigma_D \cup \Sigma_N
\\ 
\label{EQ10}  
\bu(0) = \bu_0, \ \ \bu'(0) = \bw_0, \ \ \theta(0) = \theta_0 & \mbox{\rm in} & \Omega.
\end{eqnarray}
\end{Problem}

We briefly comment on Problem~\ref{p222}. 
Equation (\ref{EQ1}) is the equation of motion while (\ref{EQ2}) is the constitutive law for viscoelastic materials in which 
${\mathscr A}$ is a nonlinear 
viscosity operator, ${\mathscr B}$ 
represents a nonlinear elasticity operator, 
${\mathscr C}$ is a linear relaxation tensor, 
and ${\mathscr{C}}_e$ is a nonlinear thermal expansion operator.
Relation (\ref{EQ3}) is the equation of heat transfer with the thermal conductivity operator ${\mathscr K}$ and a nonlinear operator ${\mathscr N}$ 
depending on the velocity. 
We have the clamped boundary condition (\ref{EQ4}) on $\Gamma_D$ and the surface
traction boundary condition (\ref{EQ5}) on $\Gamma_N$. 
Relation (\ref{EQ6}) is the multivalued contact condition with normal damped response 
in which $\partial j_\nu$ denotes the Clarke subdifferential of a given function $j_\nu$,  
$k$ is a damper coefficient, and ${\mathscr R}$ 
is a regularization operator. 
Conditions (\ref{EQ7}), (\ref{EQ8}) 
represent a version of Coulomb's law of 
dry friction in which $F_b$ 
is a friction bound.
In contrast to the conditions used in the literature, both the damper coefficient in (\ref{EQ6}) and the friction bound in 
(\ref{EQ7}), (\ref{EQ8}) depend on the temperature.
In addition, the friction bound may depend on the total slip represented by the quantity
$$
S(\bx,t)=\int_0^t\|\bu_\tau(\bx,s)\|\,ds.
$$
\noindent 
For simplicity, in (\ref{EQ9}), we assume that the temperature vanishes on $\Sigma_D \cup \Sigma_N$. Finally, initial
displacement, velocity and temperature are specified in (\ref{EQ10}). 

To state the weak formulation of Problem~\ref{p222}, we introduce the spaces
$$
\displaystyle 
V = \left\{ v \in H^1(\Omega; \real^d) 
\mid v = 0 \ \ {\rm on} \ \Gamma_D \right\},  
\ \ \  
E = \left\{ \theta \in H^1(\Omega) \mid \theta = 0 \ \ {\rm on} \ \Gamma_D \cup \Gamma_N \right\} ,
$$ 
and ${\cal H} = L^2(\Omega;\mathbb{S}^d)$. 
We use the standard notation for the Sobolev spaces, see e.g.~\cite{DMP2,MOSBOOK,Zeidler} 
and for the Bochner spaces as in Section~\ref{inequality}.
The trace on the boundary $\partial\Omega$ of an element $\bv\in H^1(\Omega; \real^d)$ is denoted by the same symbol.
We use the notation $v_\nu$ and $\bv_\tau$ for its normal and tangential traces, and 
$\| \gamma\|$ denotes the norm of the trace operator $\gamma \colon V \to L^2(\Gamma, \real^d)$.

%

\medskip

The following assumptions on the data of Problem~\ref{p222} will be needed 
throughout this section. We assume that the viscosity operator ${\mathscr A}$, 
the elasticity operator ${\mathscr B}$, the relaxation operator ${\mathscr C}$ 
and the thermal expansion operator ${\mathscr C}_e$ satisfy the following hypotheses.

\medskip

\noindent 
$\underline{H({\mathscr A})}:$ \quad 
$\displaystyle {\mathscr A} \colon Q_T \times \real \times \mathbb{S}^{d} \to \mathbb{S}^{d}$ is such that 

\lista{
	\item[(a)]  
	${\mathscr A}(\cdot, \cdot, \theta, \bvarepsilon)$ is measurable on $Q_T$ for all $\theta \in \real$, 
	$\bvarepsilon \in \mathbb{S}^{d}$. 
	\smallskip
	\item[(b)]
	${\mathscr A}(\bx,t,\cdot,\bvarepsilon)$ is $L_{\mathscr A}$--Lipschitz continuous for all 
	$\bvarepsilon \in \mathbb{S}^{d}$, 
	a.e.\ $(\bx, t)\in Q_T$ with $L_{\mathcal A} > 0$. \smallskip
	\item[(c)]
	$\| {\mathscr A} (\bx, t, \theta, \bvarepsilon) \|_{\mathbb{S}^{d}} \le a_0(\bx,t) + a_1 | \theta | + a_2 \| \bvarepsilon\|_{\mathbb{S}^{d}}$ 
	for all $\theta \in \real$, $\bvarepsilon \in \mathbb{S}^{d}$, 
	a.e.\ $(\bx, t)\in Q_T$ with $a_0 \in L^2(Q_T)$, 
	$a_0$, $a_1$, $a_2 \ge 0$. \smallskip
	\item[(d)] 
	${\mathscr A}(\bx,t, \theta,\cdot)$ is continuous and strongly monotone with $m_{{\mathscr A}} > 0$,
	for all $\theta \in \real$, a.e.\ $(\bx, t)\in Q_T$. 
	\item[(e)]
	${\mathscr A} (\bx, t, \theta, \bzero) = \bzero$  
	for all $\theta \in \real$, a.e.\ $(\bx, t)\in Q_T$.
}

\medskip

\noindent
$\underline{H({\mathscr B})}:$ \quad 
$\displaystyle {\mathscr B} \colon Q_T \times \mathbb{S}^{d} \to \mathbb{S}^{d}$ is such that 


\lista{
	\item[(a)]
	${\mathscr B}(\cdot,\cdot, \bvarepsilon)$ is measurable on $Q_T$ for all $\bvarepsilon \in \mathbb{S}^{d}$.
	\smallskip
	\item[(b)]
	$\| {\mathscr B}(\bx,t,\bvarepsilon) \|_{\mathbb{S}^{d}} \le
	b_0(\bx,t) + b_1 \, \| \bvarepsilon \|_{\mathbb{S}^{d}}$ 
	for $\bvarepsilon \in \mathbb{S}^{d}$,
	a.e. $(\bx, t) \in Q_T, b_0 \in L^2(Q_T), b_0, b_1 \ge 0$.
	\smallskip
	\item[(c)]
	${\mathscr B}(\bx,t,\cdot)$ is
	$L_{\mathscr B}$--Lipschitz continuous for 
	a.e. $(\bx, t) \in Q_T$ with $L_{\mathcal B} > 0$. 
}

\medskip

\noindent $\underline{H({\mathscr C})}:$ \quad 
${\mathscr C} \colon Q_T \times \mathbb{S}^{d} \to \mathbb{S}^{d}$ is such that 

\lista{
	\item[(a)]
	${\mathscr C}(\bx,t,\bvarepsilon) = 
	c(\bx,t) \, \bvarepsilon$ 
	for all $\bvarepsilon \in \mathbb{S}^{d}$, 
	a.e. $(\bx, t) \in Q_T$.
	\smallskip
	\item[(b)]
	$c(\bx,t) = (c_{ijkl}(\bx,t))$ with $c_{ijkl} = c_{jikl} = c_{lkij} \in L^2(0, T; L^{\infty}(\Omega))$.
}

\medskip

\noindent
$\underline{H({\mathscr C}_e)}:$ \quad 
$\displaystyle {\mathscr C}_e \colon Q_T \times \real \to \mathbb{S}^{d}$ is such that 
\smallskip

\lista{
	\item[(a)]
	${\mathscr C}_e(\cdot,\cdot, r)$ is measurable on $Q_T$ for all $r \in \real$.
	\smallskip
	\item[(b)]
	$\displaystyle 
	\| {\mathscr C}_e(\bx, t, r) \|_{\mathbb{S}^{d}} \le c_{0e}(\bx,t) + c_{1e} \, | r |$ 
	for $r \in \real$, a.e. $(\bx, t) \in Q_T,  c_{0e} \in L^2(Q_T), c_{0e}, c_{1e} \ge 0$.
	\smallskip
	\item[(c)]
	$\displaystyle 
	{\mathscr C}_e (\bx, t, \cdot)$
	is $L_e$--Lipschitz continuous for a.e. $(\bx, t) \in Q_T$ with $L_e > 0$. 
}

\medskip

The thermal conductivity operator ${\mathscr K}$, the operator ${\mathscr N}$ in the heat equation, 
and the tangential function $h_\tau$ satisfy the following assumptions.

\medskip

\noindent
$\underline{H({\mathscr K})}:$ \quad 
$\displaystyle {\mathscr K} \colon Q_T \times \real^d \to \real^d$ is such that 

\lista{
	\item[(a)]
	${\mathscr K}(\cdot,\cdot, \bxi)$ is measurable on $Q_T$ for all $\bxi \in \real^d$.
	\smallskip
	\item[(b)]
	$\displaystyle 
	\| {\mathscr K}(\bx, t, \bxi) \| \le 
	{\bar k}_0 
	(\bx,t) + {\bar k}_1 \, \| \bxi \|$ for all 
	$\bxi \in \real^d$, a.e. $(\bx, t) \in Q_T$ with  ${\bar k}_0 \in L^2(Q_T)$, 
	${\bar k}_0 \ge 0$, ${\bar k}_1 > 0$.
	\smallskip
	\item[(c)]
	${\mathscr K}(\bx, t, \cdot)$ 
	is continuous and strongly monotone with $m_{\mathscr K}>0$ 
	for a.e. $(\bx, t) \in Q_T$.  
}

\medskip

\noindent
$\underline{H({\mathscr N})}:$ \quad 
$\displaystyle {\mathscr N} \colon Q_T \times V \to L^2(\Omega)$ is such that 

\lista{
	\item[(a)]
	$\displaystyle 
	{\mathscr N}( \bx, t, \cdot)$ is 
	$L_{\mathscr N}$--Lipschitz continuous for 
	a.e. $(\bx, t) \in Q_T$ with $L_{\mathscr N} > 0$.
	\smallskip 
	\item[(b)]
	${\mathscr N}(\cdot, \cdot, v) \in L^2(Q_T)$ for all $v \in E$.
}

\medskip

\noindent
$\underline{H(h_\tau)}:$ \quad 
$\displaystyle h_\tau \colon \Gamma_C \times \real_+ \to \real_+$ is such that 

\lista{
	\item[(a)]
	$h_\tau (\bx, \cdot)$ is $L_\tau$--Lipschitz continuous for a.e. $\bx \in \Gamma_C$ with $L_\tau > 0$. \smallskip
	\item[(b)]
	$h_\tau(\cdot, r) \in L^2(\Gamma_C)$ for all $r \in \real_+$.	
}

\medskip

The potential functions $j_\nu$ and $j$, 
the damper coefficient $k$, 
the friction bound $F_b$, and the 
regularization operator ${\mathscr R}$ 
satisfy the following hypotheses.

\medskip

\noindent
$\underline{H(j_\nu)}:$ \quad 
$\displaystyle j_\nu \colon \Gamma_C \times \real 
\to \real$ is such that 

\lista{
\item[(a)]
$j_\nu (\cdot, r)$ is measurable on $\Gamma_C$ 
for all $r \in \real$ and there is 
${\bar{e}} \in L^2(\Gamma_C)$ such that 
$j_\nu(\cdot, {\bar{e}} (\cdot)) \in
L^1(\Gamma_C)$. \smallskip
\item[(b)] 
$j_\nu(\bx, \cdot)$ is locally Lipschitz on 
$\real$ for a.e. $\bx \in \Gamma_C$. 
\smallskip
\item[(c)] 
$| \partial j_\nu(\bx, r) | \le {\bar{c}}_0$
for all $r\in \real$, a.e. $\bx \in \Gamma_C$, 
with ${\bar{c}}_0\ge 0$. \smallskip
\item[(d)]
$j_\nu^0(\bx, r_1; r_2 - r_1) + j_\nu^0(\bx, r_2; r_1 - r_2) \le
{\bar{\beta}} \, | r_1 - r_2 |^2$ 
for all $r_1$, $r_2 \in \real$, a.e. 
$\bx \in \Gamma_C$ with ${\bar{\beta}} \ge 0$.
}

\medskip

\noindent
$\underline{H(j)}:$ \quad 
$\displaystyle j \colon \Sigma_C \times \real 
\to \real$ is such that 

\lista{
	\item[(a)]
	$j (\cdot, \cdot, r)$ is measurable for all 
	$r \in \real$ and there is 
	${\bar{e}}_2 \in L^2(\Gamma_C)$ such that 
	$j(\cdot, \cdot, {\bar{e}}_2 (\cdot)) \in
	L^1(\Sigma_C)$. \smallskip
	\item[(b)] 
	$j(\bx, t, \cdot)$ is locally Lipschitz on 
	$\real$ for a.e. $(\bx, t) \in \Sigma_C$. 
	\smallskip
	\item[(c)] 
	$| \partial j(\bx, t, r) | \le c_0(\bx, t) 
	+ c_1 |r|$
	for all $r\in \real$, a.e. $(\bx,t) \in \Sigma_C$,  
	$c_0 \in L^2(\Sigma_C)$, 
	$c_0, c_1 \ge 0$. \smallskip
	\item[(d)]
	$j^0(\bx, t, r_1; r_2 - r_1) + 
	j^0(\bx, t, r_2; r_1 - r_2) \le
	m_0 \, | r_1 - r_2 |^2$ 
	for all $r_1$, $r_2 \in \real$, a.e. 
	$(\bx,t) \in \Sigma_C$ with $m_0 \ge 0$.
}

\medskip

\noindent
$\underline{H(k)}:$ \quad 
$\displaystyle k \colon \Sigma_C \times \real^2 \to \real_+$ is such that 

\lista{
	\item[(a)]
	$k(\cdot, \cdot, \theta, r)$ 
	is measurable for all $\theta$, $r \in \real$.
\smallskip
	\item[(b)] 
	there exist $k_1$, $k_2$ such that 
	$0 < k_1 \le k(\bx, t, \theta, r) \le k_2$ 
	for all $\theta$, $r \in \real$, a.e. 
	$(\bx,t) \in \Sigma_C$.
	\smallskip
	\item[(c)] 
	$k(\bx, t, \cdot, \cdot)$ is $L_k$--Lipschitz continuous for a.e. $(\bx,t) \in \Sigma_C$ 
	with $L_k > 0$.
}

\medskip

\noindent
$\underline{H(F_b)}:$ \quad 
$\displaystyle F_b \colon \Sigma_C \times \real^2 \to \real$ is such that 

\lista{
	\item[(a)]
	$F_b(\cdot, \cdot, \theta, y)$ 
	is measurable for all $\theta$, $y \in \real$.
	\smallskip
	\item[(b)] 
	$F_b(\bx, t, \cdot, \cdot)$ is $L_{F_b}$--Lipschitz continuous for a.e. 
	$(\bx,t) \in \Sigma_C$ with $L_{F_b} > 0$.
	\smallskip
	\item[(c)]
	$F_b (\bx, t, 0, 0)$ belongs to $L^2(\Sigma_C)$.
}

\medskip

\noindent
$\underline{H({\mathscr R})}:$ \quad 
$\displaystyle {\mathscr R} \colon H^{-1/2}(\Gamma) \to L^2(\Gamma)$ is a linear and bounded operator. 

\medskip

For the densities of body forces, surface tractions and heat sources, and the initial data, we need the hypothesis. 

\medskip

\noindent
$\underline{(H_8)}:$ \ 
$\fb_0 \in L^2(0,T;L^2(\Omega; \real^d))$,  
$\fb_N\in L^2(0,T;L^2(\Gamma_N; \real^d))$,   
$h_0 \in L^2(0, T; E^*)$,

\qquad $\bu_0$, $\bw_0 \in V$, $\theta_0 \in E$.

\medskip

Moreover, we define $\hb_1 \in L^2(0, T; V^*)$ by
\begin{equation}\label{h1}
\langle \hb_1(t), \bv \rangle_{V^* \times V}
= 
\langle \fb_0(t), \bv \rangle_{H} +
\langle \fb_N(t), \bv \rangle_{L^2(\Gamma_N;\real^d)} 
\end{equation}
for $\bv \in V$, $t \in (0, T)$.

Note that the regularization operator 
${\mathscr R}$ in $H(\mathscr R)$ has been introduced and used earlier in several papers on contact probems, see~\cite{Denk2006,DUV,MG2019,SST}.
Following a standard procedure, see e.g.~\cite{MOSBOOK,MIGSZA,SOFMIG}, we obtain 
the weak formulation of Problem~\ref{p222}.    
\begin{Problem}\label{p333}
Find a displacement field $\bu \colon (0, T) \to V$ and a temperature $\theta \colon (0, T) \to E$ such that for a.e.\ $t \in (0, T)$,
	\begin{eqnarray}
	&&
	\langle \bu''(t), \bv - \bu'(t) \rangle_{V^*\times V} + 
	({\mathscr A}(t, \theta(t),\bvarepsilon({\bu}'(t))),
	\bvarepsilon(\bv-{\bu}'(t)))_{\cal H} \nonumber \\[1mm]
	&&\ \ {}+\Big({\mathscr B}(t,\bvarepsilon({\bu}(t)))  
	+\int_0^t{\mathscr C}(t-s)
	\bvarepsilon({\bu}'(s))\,ds 
	+ {\mathscr C}_e(t,\theta(t)),
	\bvarepsilon(\bv-{\bu}'(t))\Big)_{\cal H} \nonumber \\[1mm]
	&&\ \ {} +\int_{\Gamma_C}F_b\Big(t, {\mathscr R} \theta(t),
	\int_0^t\|\bu_\tau(s)\|\,ds \Big) \,
	\big(\|\bv_\tau\|-\|\bu_\tau'(t)\|\big) \, d\Gamma\nonumber\\[1mm]
	&&\ \  {} +\int_{\Gamma_C} 
	k(t, {\mathscr R}\theta(t), u_\nu(t))\, j_\nu^0 ({u}_\nu'(t); v_\nu-{u}_\nu'(t))\, d\Gamma
	\ge \langle\hb_1(t),\bv-{\bu}'(t)\rangle_{V^*\times V} \label{P}
	\end{eqnarray}
for all $\bv \in V$, 
and
\begin{eqnarray}
&&\langle \theta'(t), \zeta \rangle_{E^* \times E}  
+\langle 
{\mathscr K}(t, \nabla \theta(t)), \nabla \zeta \rangle_{H}
+\int_{\Gamma_C}
j^0(t, \theta(t); \zeta)\, d\Gamma \nonumber
\\
&&
\ \ {} 
\ge \langle 
{\mathscr N}(t, \bu'(t)) + h_0(t), \zeta
\rangle_{E^* \times E}
+ \int_{\Gamma_C}
h_\tau(\| \bu'_{\tau}\|) \zeta \, d\Gamma
\label{p334}
\end{eqnarray}
for all $\zeta \in E$
with	
\begin{equation}
\bu(0)=\bu_0,\ \ 
{\bu}'(0)=\bw_0,\ \  
\theta(0) = \theta_0.\label{Initial}
\end{equation}
\end{Problem}

\smallskip

The unique solvability of Problem~\ref{p333} is provided in the following result.

\begin{Theorem}\label{t333}
Assume the hypotheses
${H(\mathscr A)}$, ${H(\mathscr B)}$, 
${H(\mathscr C)}$, ${H(\mathscr C_e)}$, 
${H(\mathscr K)}$, ${H(\mathscr N)}$, $H(h_\tau)$,
$H(j_\nu)$, $H(j)$, $H(k)$, $H(F_b)$, 
${H(\mathscr R)}$, and $(H_8)$.
If 
\begin{equation}\label{sss}
m_{\mathscr A} > {\bar{\beta}} \, k_2 \| \gamma \|^2
\ \ \mbox{and} \ \ 
m_{\mathscr K} > m_0 \|\gamma\|^2 ,
\end{equation}
then Problem~$\ref{p333}$ has a unique solution with regularity $\bu \in {\mathcal V}$, 
$\bu' \in {\mathcal W}$, and 
$\theta \in {\mathbb E}$.
\end{Theorem}

\noindent
{\bf Proof.}  
We formulate Problem~$\ref{p333}$ in the form of Problem~\ref{SYSTEM} and apply Theorem~\ref{MAIN}. To this end, we introduce the following notation. 

We work with two evolution triples of spaces 
$(V, H, V^*)$ with $H = L^2(\Omega;\real^d)$, 
and $(E, X, E^*)$ with $X = L^2(\Omega)$. 
Let $Y=Z=L^2(\Gamma_C)$. 
We introduce the operators 
$A \colon (0,T)\times X \times V \to V^*$, 
$R \colon {\mathcal V} \to L^2(0, T; Y)$, 
$R_1 \colon {\mathcal V} \to {\mathcal V^*}$, 
and 
$S \colon {\mathcal V} \to L^2(0, T; Z)$ 
defined by
\begin{eqnarray}
&&\hspace{-9mm}\label{b1} \langle 
A (t, \theta,\bv), \bzeta \rangle_{V^*\times V} =
({\mathscr A} (t, \theta, \bvarepsilon(\bv)) 
+ {\mathscr C}_e(t,\theta), \varepsilon(\bzeta) \rangle_{\mathcal H} \label{operA}\\[1mm]
&&
\qquad \qquad 
\qquad \mbox{for}\ \theta \in X, \, \bv,\bzeta \in V,  \, t\in (0,T), \nonumber \\ [1mm] 
&&\hspace{-9mm}\label{b3} 
(R\bv)(t) = 
\int_0^t\Big\|\int_0^s \bv_\tau(r)\,dr+\bu_{0\tau}\Big\| \,ds
\quad \mbox{for} \ \bv \in{\cal V},\  t \in (0,T), \\[1mm]
&&\hspace{-9mm} \langle 
(R_1 \bv)(t), \bzeta \rangle_{V^*\times V} = 
\Big( {\mathscr B}\big(t, 
\big(\int_0^t\bvarepsilon({\bv}(s))\,ds+\bu_0\big)\big),
\bvarepsilon(\bzeta) \Big)_{\cal H}
\nonumber \\[1mm]
&&\hspace{-9mm}\ \ {} +\Big(\int_0^t {\mathscr C} (t-s)\bvarepsilon({\bv}(s))\,ds,
\bvarepsilon(\bzeta)\Big)_{\cal H} 
\ \ 
\mbox{for} \ \bv \in{\cal V},\ \bzeta \in V, \ t \in (0, T),\label{b2}\\[1mm]
&&\hspace{-9mm}\label{b4} 
(S\bv)(t) = \int_0^t v_\nu(s)\,ds+u_{0\nu} 
\quad \mbox{for} \ \bv \in{\cal V},\  t \in (0,T).
\end{eqnarray} 
Further, we define functions 
$J \colon (0, T) \times X \times Z \times V \to\real$ and 
$\varphi \colon (0, T) \times X \times Y \times V \to \real$ by
\begin{eqnarray}
&&\label{b5} 
J(t,\theta,z,\bv)=\int_{\Gamma_C} 
k(t,{\mathscr R} \theta, z)\, j_\nu (v_\nu) \,d\Gamma\quad\mbox{for}
\ \theta \in X, z\in Z,\ \bv\in V,\ t\in (0, T),\\ [1mm]
&&\label{b6} 
\varphi(t,\theta, y,\bv)=\int_{\Gamma_C} 
F_b(t, {\mathscr R} \theta,y) \, \| \bv_\tau \| \,d\Gamma\quad\mbox{for}\ \theta \in X, y\in Y,\ 
\bv\in V,\ t \in (0,T),
\end{eqnarray}
the operator 
$B \colon (0, T) \times V \times E \to E^*$
and the function 
$g \colon (0, T) \times V \times E \to \real$ 
by
\begin{eqnarray}
&&\label{b52} 
\langle B(t, \bv, \theta), \zeta \rangle_{E^* \times E} 
= \langle {\mathscr K}(t, \nabla \theta), 
\nabla \zeta \rangle_H - 
\langle {\mathscr N}(t,\bv), 
\zeta \rangle_{E^* \times E} 
- \langle h_\tau (\| \bv_\tau \|), \zeta 
\rangle_{L^2(\Gamma_C)}
\\[1mm]
&&
\qquad\qquad\qquad\ \ \ 
\quad\mbox{for} \ \bv\in V, \,  
\theta, \, \zeta \in E, \, t\in (0, T),
\nonumber \\[1mm]
&&\label{b62} 
g (t, \bv, \theta)=\int_{\Gamma_C} 
j(t, \theta) \, d\Gamma
\quad\mbox{for}\ \theta \in E, \ t \in (0,T).
\end{eqnarray}
\noindent 
Observe that the operator $B$ is independent 
on the history, and $g$ does not depend on $\bv$.
We denote 
$\bw = \bu'$, i.e., 
\begin{equation}\label{solutionu}
\bu(t)=\int_0^t\bw(s)\,ds+\bu_0\ \ \mbox{for all} \ t\in [0,T].
\end{equation}
With the notation above we consider the following problem in terms of the velocity 
and the temperature.
\begin{Problem}\label{Problem2n} 
	Find $\bw \in {\mathbb{W}}$ and $\theta \in {\mathbb E}$ such that
	\begin{equation*}
	\begin{cases}
	\displaystyle
	\bw'(t) + A(t, \theta(t), \bw(t)) + 
	({R}_1 \bw)(t) + 
	\partial J (t, \theta(t), ({\mathcal S} \bw)(t), \bw(t)) \\[1mm]
	\qquad\qquad\qquad\qquad \quad \ \   
	+ \, \partial_c \varphi (t, \theta(t), 
	({R} \bw)(t), \bw(t)) \ni \hb_1(t) 
	\ \ \mbox{\rm a.e.} \ \ t \in (0,T), \\[1mm]
	\theta'(t) + B(t, \bw(t), \theta(t))
	+ \partial g(t, \bw(t), \theta(t)) \ni h_0(t) \ \ \mbox{\rm a.e.} \ \ t \in (0,T), \\[1mm]
	\bw(0) = \bw_0, \ \theta(0) = \theta_0 .
	\end{cases}
	\end{equation*}
\end{Problem}

We apply Theorem~\ref{MAIN} in studying Problem~\ref{Problem2n}. For this purpose, we check that the hypotheses 
$H(A)$, $H(J)$, $H(\varphi)$, $H(B)$, $H(g)$, $(H_3)$--$(H_5)$ are satisfied. 

First, we note that by hypotheses $H({\mathscr A})$(a), (c), and $H({\mathscr C}_e)$(a), (b), 
the operator $A$ defined by (\ref{operA}) 
satisfies conditions $H(A)$(a) and (c).
It follows from the proof of~\cite[Theorem~7.3]{MOSBOOK} that 
$A(t, \theta, \cdot)$ is demicontinuous for $\theta
\in X$, a.e. $t \in (0, T)$.
Next, from $H({\mathscr A})$(b) and 
$H({\mathscr C}_e)$(c), we deduce that 
$A(t, \cdot, \bv)$ is 
$(L_{\mathscr A} + L_e)$--Lipschitz for all $\bv \in V$, a.e. $t \in (0, T)$. 
Moreover, since 
$A(t, \theta, \cdot)$ is strongly monotone
with $m_{\mathscr A}>0$ for $\theta \in X$, 
a.e. $t \in (0, T)$, by Remark~\ref{Remark22}, 
we infer that $H(A)$(d) holds.
We conclude that the operator $A$ satisfies 
$H(A)$ with $m_A = m_{\mathscr A}$.

Exploiting the hypotheses $H(k)$, $H(j_\nu)$ and 
$H(\mathscr R)$, after some calculation, 
we can obtain, similarly as in the proof of~\cite[Theorem~13]{HMS2017}, that
the function $J$ defined by (\ref{b5}) satisfies $H(J)$ with $m_J = 
{\bar{\beta}} \, k_2 \| \gamma \|^2$.
Furthermore, the function $\varphi$ defined by (\ref{b6}) under the hypotheses $H(F_b)$ and $H(\mathscr R)$ satisfy $H(\varphi)$. 
For details, we refer to~\cite[Theorem~13]{HMS2017}.

It follows from $H(\mathscr K)$(a) and 
$H(\mathscr N)$(b) that $H(B)$(a) holds. 
The hypotheses 
$H(\mathscr N)$(a) and $H(h_\tau)$(a) imply that 
$B(t,\cdot,\theta)$ is continuous for all 
$\theta \in E$, a.e. $t \in (0,T)$. 
Hence $H(B)$(b) is satisfied. 
The growth condition $H(B)$(c) is a consequence 
of
$H({\mathscr K})$(b), $H(h_\tau)$, and 
$H({\mathscr N})$.
By~\cite[Lemma~6]{MIGSZA}, we know that 
$B(t,\bv, \cdot)$ is pseudomonotone for all $\bv \in V$, a.e. $t \in (0, T)$. So,  from~\cite[Theorem~3.69(ii)]{MOSBOOK}, 
we know that $B(t,\bv, \cdot)$ is demicontinuous.
Further, using $H(\mathscr N)$ and $H(h_\tau)$, 
we infer that 
$B(t,\cdot, \theta)$ is Lipschitz for all 
$\theta \in E$, a.e. $t \in (0,T)$.
Combining the latter with the fact that 
${\mathscr K}(\bx, t, \cdot)$ is strongly monotone 
for a.e. $(\bx, t) \in Q_T$, by Remark~\ref{Remark22}, 
we deduce that $H(B)$(d) holds with 
$m_B = m_{{\mathscr K}}$. Thus $H(B)$ is verified.

The conditions $H(g)$(a) and (b) hold trivially. 
Using $H(j)$(a)--(c), it follows from~\cite[Theorem~3.47(iii)]{MOSBOOK} that 
$g(t, \bv, \cdot)$ is Lipschitz on bounded sets for a.e. $t \in (0, T)$ which implies $H(g)$(c).
By $H(j)$(c), we obtain $H(g)$(d) analogously as in~\cite[Theorem~3.47(vi)]{MOSBOOK}. 
Next, using $H(j)$(d), we have 
$$
g^0(t, \bv_1, \theta_1; \theta_2 - \theta_1) 
+ g^0(t, \bv_2, \theta_2; \theta_1 - \theta_2) 
\le m_g \, \| \theta_1 - \theta_2 \|^2_E
$$
for all $\theta_1$, $\theta_2 \in E$, 
a.e. $t \in (0, T)$ with $m_g = m_0 \|\gamma\|^2$. 
This entails the condition $H(g)$(e). 
Hence $H(g)$ holds.

Condition $(H_3)$ follows from (\ref{h1}) and $(H_8)$, and the smallness condition 
$(H_4)$ is a consequence of the inequalities (\ref{sss}). 
The operators 
$R$, $R_1$ and $S$ defined by (\ref{b3}), (\ref{b2}) 
and (\ref{b4}) are history-dependent, i.e.,
they satisfy conditions $(H_5)$(a), (b) and (d), respectively. These conditions can be checked analogously as in~\cite{AMMA14}. 

We have verified all hypotheses of Theorem~\ref{MAIN}. 
Then, we deduce that Problem~\ref{Problem2n} has a unique solution $(\bw, \theta) \in {\mathbb W} \times {\mathbb E}$.
Next, we recall that the function 
$\bu \colon(0,T) \to V$ is defined 
by (\ref{solutionu}). 
By a direct calculation we can check that Problems~$\ref{p333}$ and~\ref{Problem2n}
have unique solutions. This entails that 
$(\bw, \theta) \in {\mathbb W} 
\times {\mathbb E}$ solves Problem~\ref{Problem2n} if and only if 
$\bu \in {\mathcal V}$, 
$\bu' \in {\mathcal W}$, and 
$\theta \in {\mathbb E}$ solves Problem~\ref{p333}.
This completes the proof of the theorem.
{\hfill $\square$}

\medskip

A triple of functions $(\bu, \theta, \bsigma)$ which satisfies (\ref{P}), (\ref{p334}), (\ref{Initial}) and (\ref{EQ2})
is called a {\it weak solution} to Problem~\ref{p222}.
We conclude that, under the assumptions of Theorem~\ref{t333},
Problem~\ref{p222} has a unique weak solution. Moreover, the solution has the regularity
$$
\bu \in  L^2(0,T;V), \quad 
\bu' \in {\mathbb W},\quad
\theta \in {\mathbb E}, \quad
\bsigma \in L^2(0, T; {\cal H}), 
\quad {\rm Div} \, \bsigma \in L^2(0,T;V^*).
$$

Finally, we note that exploting Theorem~\ref{MAIN}, 
new results on the unique solvability of other contact problems can be obtained. 
We refer to~\cite{MMM,SMH2018} and \cite[Section~7.1]{MOSBOOK} for examples of nonconvex superpotentials which satisfy hypotheses $H(j_\nu)$ and $H(j)$.

\end{document}